\documentclass[11pt]{article}
 
\usepackage{varioref,graphicx, epstopdf, amsmath,subfigure, amsfonts} 

\usepackage{color,hyperref} 

\catcode`\@=11
\@addtoreset{equation}{section}
 \catcode`\@=12

\newtheorem{Theorem}{Theorem}[section]

\newtheorem{Proposition}{Proposition}[section]
\newtheorem{Lemma}{Lemma}[section]
\newtheorem{Remark}{Remark}[section]
\newtheorem{Corollary}{Corollary}[section]

\def\dd{{\mathsf d}}
\def\MD{{\mathsf D}}
\def\dist{{\mathsf{d}}}
\def\diam{{\mathsf{diam}}}

\def\p{\lambda_{0}}
\def\a{{\bf a}}

\def\MC{{\mathsf C}}

\def\Proof{\noindent{{\bf Proof. }}}
\def\square{\vbox{ \hrule height .4pt
\hbox{\vrule width .4pt height 7pt \kern 7pt \vrule width
.4pt}\hrule height .4pt }}
\def\QED{\hfill {$\square$}\bigskip}

\def\fa{\forall}
\def\hb{\hfill\break}

\def\norm#1{\| #1\|}

\def\tfrac#1#2{{\textstyle\frac{#1}{#2}}}

\def\d{{\delta}}

\def\G{{\mathcal V}}

\def\f{\varphi_{c}}

\def\ow{{\overline w}}

\def\uw{{\underline w}}

\def\CC{{\cal{C}}}
\def\MM{{\cal{M}}}
\def\KK{{\cal{K}}}
\def\HH{{\cal{H}}}

\def\SS{{\cal{H}}}

\def\R{{\mathbb{R}}}
\def\N{{\mathbb{N}}}
\def\Z{{\mathbb{Z}}}

\title{Stationary layered  solutions for a system of Allen-Cahn type equations }
\author{\small
Francesca Alessio}

\date{ }

\begin{document}
\small

\maketitle

\centerline{
{\scriptsize
{\baselineskip=8pt{\footnote{\scriptsize{Supported by MURST 
Project `Metodi Variazionali ed 
Equazioni Differenziali Non Lineari'}}}}
 Dipartimento di Ingegneria Industriale e Scienze
Matematiche, Universit\`a Politecnica delle Marche}
}

\centerline{
{\scriptsize Via
Brecce Bianche, I--60131 Ancona, Italy }}

\centerline{
{\scriptsize
{\sf e-mail}  alessio@dipmat.univpm.it}
}
\vskip4truecm

\noindent{
\small{\bf Abstract.}
We consider a class of
semilinear elliptic system of the form
\begin{equation}\label{eq:abs}
	-\Delta u(x,y)+\nabla W(u(x,y))=0,\quad (x,y)\in\R^{2},
\end{equation}
where 
$W:\R^{2}\to\R$ is a double well nonnegative symmetric potential.  We show, via variational methods, that 
if the set of solutions to the one dimensional system
$-\ddot q(x)+\nabla W(q(x))=0,\  x\in\R$, which connect the two minima of $W$ as $x\to\pm\infty$
has a discrete structure, then
(\ref{eq:abs}) has infinitely many layered solutions.}


\bigskip


\vfill\eject

\section{Introduction}
We consider semilinear elliptic system of the form
\begin{equation}\label{eq:eq0}
	-\Delta u(x,y)+\nabla W(u(x,y))=0,\qquad(x,y)\in\R^{2}
\end{equation}
where $W\in\CC^{2}(\R^{2},\R)$ satisfies
\begin{description}

\item[$(W_{1})$] there exist $\a_{\pm}\in\R^{2}$ such that $W(\a_{\pm})=0$, $W(\xi)> 0$ for every $\xi\in\R^{2}\setminus\{\a_{\pm}\}$ and $D^{2}W(\a_{\pm})$ are positive definite;

\item[$(W_{2})$] 
there exists $R>1$ such that { $\inf\{W(\xi),\, |\xi|>R\}=\mu_0>0$};

\item[$(W_{3})$] $W(-x_{1},x_{2})=W(x_{1},x_{2})$ for all $(x_{1},x_{2})\in\R^{2}$;
\end{description}
 In the sequel, without loss of generality, we will assume that ${\bf a}_{\pm}=(\pm 1,0)$.
\bigskip

The system (\ref{eq:eq0}) is the rescaled stationary system associated to the reaction-diffusion system 
\begin{equation}\label{eq:eqt}
\partial_tu(x,y)-\varepsilon^2\Delta u(x,y)+\nabla W(u(x,y))=0,\quad (x,y)\in\Omega\subset\R^2
\end{equation} 
which describes two phase physical systems or grain boundaries in alloys. As  $\varepsilon\to 0^+$, solutions to (\ref{eq:eqt}) tends almost everywhere to global minima of $W$ and sharp phase interfaces appear (see e.g. \cite{[BronsardReitich]}, \cite{[Sternberg]} and \cite{[Sternberg2]}). Then, the expansion of such solutions in a point on the interface presents, as first term, the system (\ref{eq:eq0}). From this point of view, two layered transition solutions correspond to solutions $u$ of (\ref{eq:eq0}) satisfying the asympotic conditions
\begin{equation}\label{eq:ete} \lim_{x\to\pm\infty}u(x,y)=\a_{\pm}\quad\hbox{uniformly w.r.t. }  y\in\R.
\end{equation}

S. Alama, L. Bronsard and C. Gui in \cite{[AlamaBronsardGui]} studied the existence of solution to (\ref{eq:eq0}) which satisfies the asymptotic condition (\ref{eq:ete}) for $x\to\pm\infty$ while as $y\to\pm\infty$ tends to two different one dimensional stationary waves, that is, solutions to the one dimensional associated problem 
\begin{equation}\label{eq:1dim}
		\begin{cases}-\ddot q(x)+\nabla W(q(x))=0,&x\in\R\\
			{\displaystyle\lim_{t\to\pm\infty}q(t)=\a_{\pm}}.&\end{cases}
			\end{equation}
which are furthermore minima of the action
$$
	 V(q)=\int_{\R}\tfrac12|\dot q|^{2}+W(q)\, dx
$$
over the class of trajectories connecting $\bf a_\pm$ as $x\to\pm\infty$. Such solutions are found under conditions $(W_1)$ and $(W_3)$, requiring a fast growth at infinity and assuming that there exist a finite number $k\geq 2$ of geometrically distinct one dimensional minimizing heteroclinic connections. See also the paper by Alikakos and Fusco, \cite{[AlikFusco0]} for related assumptions and existence results.

In \cite{[Schatzman]} M. Schatzman proved the same result, considering a non symmetric potential, assuming that there exists two geometrically distinct one dimensional  he\-te\-roclinic connections which are supposed to be non degenerate, i.e. the kernel of the corresponding linearized operators are one dimensional.\\
We cite also the case of a triple-well potential over $\R^2$ which was studied by L. Bronsard, C. Gui and M. Schatzman \cite{[BronGuiSchatz]} where it is proved the existence of entire solutions to (\ref{eq:eq0}), known as triple-junction solutions, which connect the three global minima of the potential $W$ in certain directions at infinity. See also \cite{[GuiSchatzman]} for quadruple junction solutions and \cite{[AlikFusco]} (and the reference therein) for potential defined over $\R^n$ with general reflection group of symmetry.

\smallskip

If $u$ is scalar valued, much is known about the corresponding heteroclinic problem (\ref{eq:eq0})-(\ref{eq:ete}). In this scalar setting, E. De Giorgi in\cite{[DeGiorgi]} has conjectured that  any entire bounded solution of $\Delta u=u^{3}-u$  with $\partial_{x_{1}}u(x)> 0$ in $\R^{n}$ for $n\le 8$
is in fact one-dimensional, i.e., modulo space rotation-traslation, it coincides   
with the unique solution 
of the one dimensional heteroclinic problem
 \begin{equation}\label{eq:ODE}
	\begin{cases}\ddot q(x)=q(x)^{3} -q(x),&
	x\in\R,\\
	 q(0)=0\hbox{ and }q(\pm\infty)=\pm 1,&
	\end{cases}
\end{equation}

 The conjecture has been proved for $n=2$ by N. Ghoussoub and C. Gui in  \cite{[GhossoubGui]} and then  by L. Ambrosio and X. Cabr\`e in \cite{[AmbrosioCabre]} (see also \cite{[AlbertiAmbrosioCabre]}) for $n=3$, even for more general double well potentials $W$. A further step in the proof of the De Giorgi conjecture has been done by O. Savin
 in \cite{[Savin]} where, for $n\leq 8$, the same one dimensional symmetry is obtained for solutions $u$ such that $\partial_{x_{1}}u(x)>0$ on $\R^{n}$ and $\lim_{x_{1}\to\pm\infty}u(x)=\pm 1$ for all $(x_{2},x_{3},...,x_{n})\in \R^{n-1}$ (see \cite{[BarlowBassGui]}, \cite{[BerestyckiHamelMonneau]}, \cite{[Farina]}, \cite{[FarinaSciunziValdinoci]} and \cite{[FarinaValdinoci1]}  for related problems). That result is completed in \cite{[DelPinoKowalczykWeiCR]}, \cite{[DelPinoKowalczykWei]} where the existence of entire solutions without any one dimensional symmetry which are increasing and asymptotic to $\pm 1$ with respect to the first variable is proved in dimension $n>8$.  \bigskip

In this paper we consider problem (\ref{eq:eq0})-(\ref{eq:ete}) and using a global variational procedure we prove that if the minimal set of one dimensional  heteroclinic connections satisfies a suitable discreteness assumption then there exist infinitely many solutions to the problem with prescribed {\it energy}.\smallskip

To explain precisely our result and to give an idea of our procedure, let us begin considering the problem already considered in \cite{[AlamaBronsardGui]} and \cite{[Schatzman]}.\\
\noindent Assuming that the minimal set of one dimensional  heteroclinic connections $\MM$ satisfies the discreteness assumption
$$
\hbox{ ($*$) $\qquad\MM=\MM^{+}\cup \MM^{-}$ with $\hbox{dist}_{L^{2}}(\MM^{+},\MM^{-})>0$.}
$$

\noindent  we will look for bidimensional solution  with prescribed different asymptots as $y\to\pm\infty$:
 \begin{equation}\label{eq:eteY}
 \hbox{dist}_{L^{2}}(u(\cdot,y),\MM^{\pm})\to 0\hbox{ as } y\to\pm\infty.
 \end{equation}

\begin{center}
\def\svgwidth{0.8\columnwidth}
 \begingroup
  \makeatletter
  \providecommand\color[2][]{%
    \errmessage{(Inkscape) Color is used for the text in Inkscape, but the package 'color.sty' is not loaded}
    \renewcommand\color[2][]{}%
  }
  \providecommand\transparent[1]{%
    \errmessage{(Inkscape) Transparency is used (non-zero) for the text in Inkscape, but the package 'transparent.sty' is not loaded}
    \renewcommand\transparent[1]{}%
  }
  \providecommand\rotatebox[2]{#2}
  \ifx\svgwidth\undefined
    \setlength{\unitlength}{482.84491367pt}
  \else
    \setlength{\unitlength}{\svgwidth}
  \fi
  \global\let\svgwidth\undefined
  \makeatother
  \begin{picture}(1,0.54258209)%
    \put(0,0){\includegraphics[width=\unitlength]{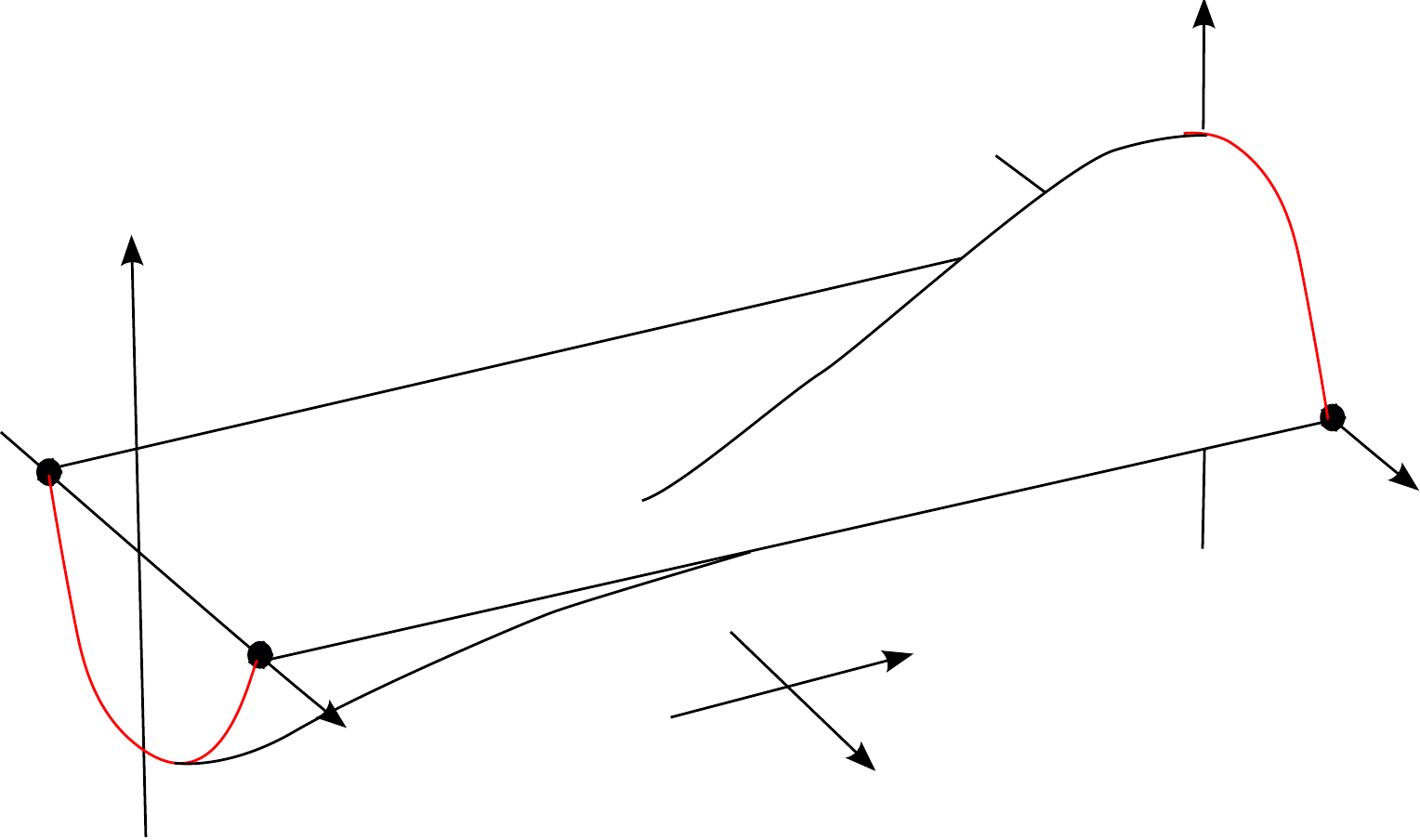}}%
    \put(-0.10165037,0.24722877){\color[rgb]{0,0,0}\makebox(0,0)[lb]{\smash{$(-1,0)$}}}%
    \put(0.11884164,0.37271384){\color[rgb]{0,0,0}\makebox(0,0)[lb]{\smash{$q_2$}}}%
    \put(0.23463616,0.0292194){\color[rgb]{0,0,0}\makebox(0,0)[lb]{\smash{$q_1$}}}%
    \put(0.01522724,0.0337552){\color[rgb]{0,0,0}\makebox(0,0)[lb]{\smash{${\cal M}^-$}}}%
    \put(0.939426934,0.43319193){\color[rgb]{0,0,0}\makebox(0,0)[lb]{\smash{${\cal M}^+$}}}%
    \put(0.87941713,0.55320585){\color[rgb]{0,0,0}\makebox(0,0)[lb]{\smash{$q_2$}}}%
    \put(0.97961809,0.200492){\color[rgb]{0,0,0}\makebox(0,0)[lb]{\smash{$q_1$}}}%
    \put(0.6256341,0.08920094){\color[rgb]{0,0,0}\makebox(0,0)[lb]{\smash{$y$}}}%
    \put(0.56505898,0.03531639){\color[rgb]{0,0,0}\makebox(0,0)[lb]{\smash{$x$}}}%
    \put(0.9725534,0.30004172){\color[rgb]{0,0,0}\makebox(0,0)[lb]{\smash{$(1,0)$}}}%
  \end{picture}%
\endgroup
 \end{center}

\noindent Note that condition $(*)$ (as the discreteness assumption made in \cite{[AlamaBronsardGui]} and \cite{[Schatzman]}) does not hold in the scalar case, where the minimal set of one dimensional solutions $\MM$  is a 
continuum homeomorphic to $\R$, being constituted by the 
translation of the unique heteroclinic solution of (\ref{eq:ODE}).\bigskip

Under our assumption, bidimensional solutions satisfying (\ref{eq:eteY}) can be found using a global variational approach (instead of the approximating procedure used in \cite{[AlamaBronsardGui]} and \cite{[Schatzman]}), considering a renormalized action functional over a suitable space. So let
$$
	 \Gamma=\{q-z_{0}\in H^{1}(\R)^{2}\,|\, \ q(x)_{1}=-q(-x)_{1},\, q(x)_{2}=q(-x)_{2}\},
	 $$
	 where $q(x)=(q(x)_{1},q(x)_{2})$ and $z_{0}$ is fixed in such a way that $z_{0}(x)_{1}=-z_{0}(-x)_{1}$,  $z_{0}(x)_{2}=z_{0}(-x)_{2}$ and $z_{0}(x)=(1,0)$ for $x>1$, be the space of one dimensional heteroclinic connections.
Setting $m=\inf_{\Gamma}V(q)$ we have that  
$$
\MM=\{q\in\Gamma\,|\, V(q)=m\}.
$$
Although there is a canonical Euler Lagrange functional associated to (\ref{eq:eq0}), 
such functional is always infinite over the solutions which we seek so, as in  \cite{[R0]} and \cite{[R1]} for
Hamiltonian ODE systems and in \cite{[AJM]} for scalar Allen-Cahn equations, we have that bidimensional heteroclinic solutions of (\ref{eq:eq0}) may be obtained as global minimizers of the ÒrenormalizedÓ  action functional
$$
\varphi(u)=\int_{\R}\frac 12\|\partial_{y}u(\cdot,y)\|^{2}_{L^{2}(\R)^{2}}+\left(V(u(\cdot,y))-m\right)\,dy
$$
which is well defined on the space
$$
{\cal H}=\{ u\in H^{1}_{loc}(\R^{2},\R^{2})\,|\, u(\cdot,y)\in\Gamma\hbox{ for almost every }y\in\R\}
$$ 
and such that $\varphi(q)=0$ for all $q\in\cal M$. We are interested on solutions which satisfies the right asympotic conditions as $y\to\pm+\infty$ that can be found as minima of $\varphi$ over the space
$$
{\cal H}_{m}=\{ u\in {\cal H}\,|\, \liminf_{y\to\pm\infty}\hbox{dist}_{L^{2}}(u(\cdot,y),\MM^{\pm})=0\}.
$$
In fact, for all $u\in\HH$ we have
$$
 \|u(\cdot,y_{1})-u(\cdot,y_{2})\|_{L^{2}(\R)^{2}}^{2}\leq (y_{2}-y_{1})\int_{y_{1}}^{y_{2}}\|\partial_{y}u(\cdot,y)\|_{L^{2}(\R)^{2}}^{2}\,dy
$$
and in particular, if $\varphi(u)<+\infty$ then the map $y\in\R\to u(\cdot,y)\in\Gamma$ is continuous with respect to the $L^{2}$ metric. Moreover,
 if $y_{1}<y_{2}$ and $u\in\HH$ then
$$
\varphi(u)\geq \left( 2\tfrac 1{y_{2}-y_{1}}\int_{y_{1}}^{y_{2}} \left( V(u(\cdot,y))-m\right)\, dy\right)^{1/2}\| u(\cdot,y_{1})-u(\cdot,y_{2})\|_{L^{2}(\R)^{2}}.
$$
By the previous estimate, if $u\in\HH_m$ we have control of the transition time from $\MM^{-}$ to $\MM^{+}$ and so concentration in the $y$ variable. Together with the symmetry in the $x$ variable this allows to get compactness of minimizing sequences and hence to prove the existence of at least one bidimensional solution to (\ref{eq:eq0}) in $\HH_{m}$, as in the Theorem by Alama Bronsard and Gui but in a slightly more general setting

\begin{Theorem}\label{T:mainPRIMO}
If $(W_{1})$-$(W_{3})$ and $(*)$ hold, then there exists $u\in{\mathcal C}^{2}(\R^{2},\R^{2})$ solution to (\ref{eq:eq0}) such that $u(x,y)\to \a_{\pm}$ as $x\to \pm\infty$ uniformly w.r.t. $y\in\R$ and
$$
\lim_{y\to\pm\infty}\dd_{H^{1}(\R)^{2}}(u(\cdot,y),\MM^{\pm})=0.
$$

\end{Theorem}

\bigskip

Note that if $u\in{\HH}$ solves the system (\ref{eq:eq0}) then 
$$
\partial^{2}_{y}u(x,y)=\underbrace{-\partial^{2}_{x}u(x,y)+\nabla W(u(x,y))}_{\displaystyle{V'(u(\cdot,y))}}
$$
In other words $u$ defines a trajectory 
$y\in\R\mapsto u(\cdot,y)\in\Gamma$ 
solution to the infinite dimensional Lagrangian system
$$
\frac {d^{2}}{dy^{2}} u(\cdot,y)=V'(u(\cdot,y))
$$
which has as equilibria point the one dimensional solution $q\in\MM$. From such point of view, bidimensional solutions in $\HH_{m}$ are  heteroclinic type solutions:

 \begin{center}
 \def\svgwidth{0.8\columnwidth}
 \begingroup
  \makeatletter
  \providecommand\color[2][]{%
    \errmessage{(Inkscape) Color is used for the text in Inkscape, but the package 'color.sty' is not loaded}
    \renewcommand\color[2][]{}%
  }
  \providecommand\transparent[1]{%
    \errmessage{(Inkscape) Transparency is used (non-zero) for the text in Inkscape, but the package 'transparent.sty' is not loaded}
    \renewcommand\transparent[1]{}%
  }
  \providecommand\rotatebox[2]{#2}
  \ifx\svgwidth\undefined
    \setlength{\unitlength}{428.3480957pt}
  \else
    \setlength{\unitlength}{\svgwidth}
  \fi
  \global\let\svgwidth\undefined
  \makeatother
  \begin{picture}(1,0.47325704)%
    \put(0,0){\includegraphics[width=\unitlength]{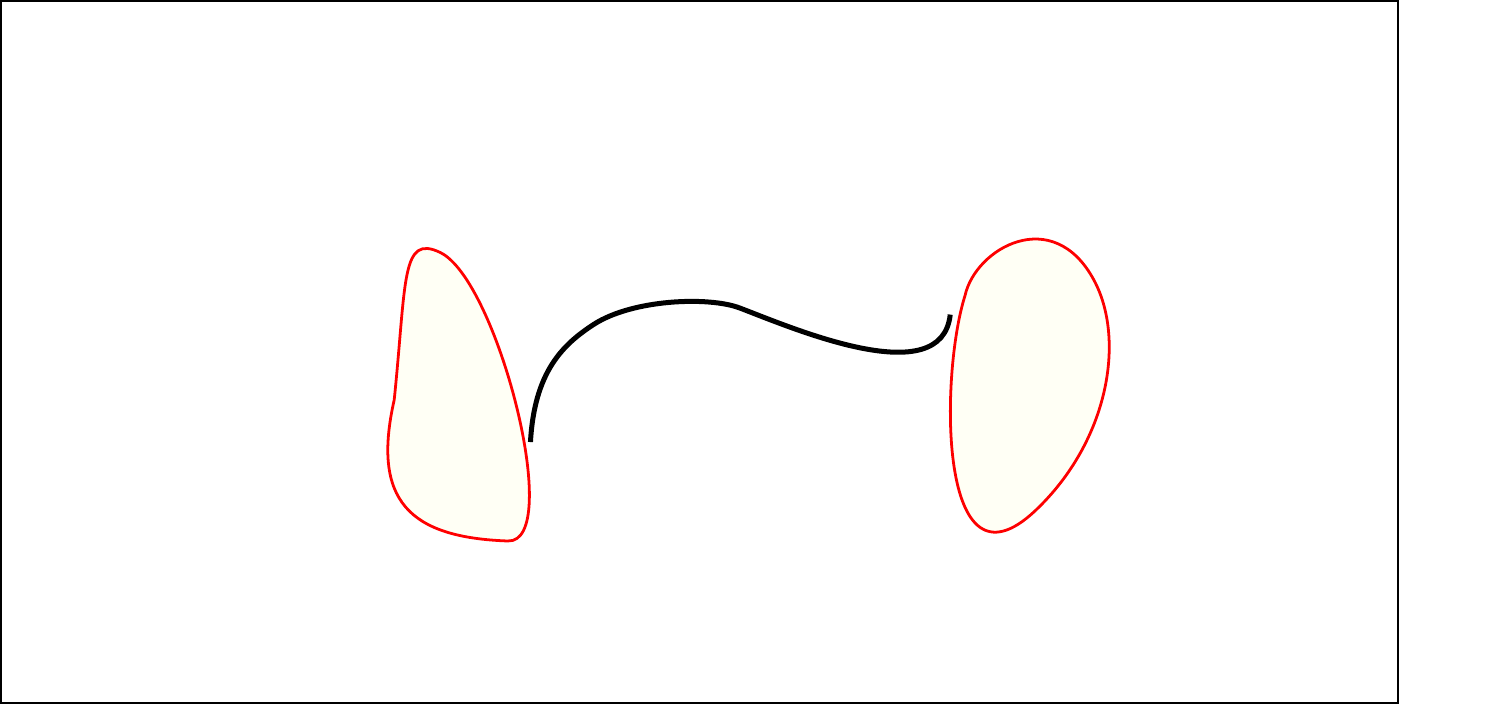}}%
    \put(0.17531519,0.20301057){\color[rgb]{0,0,0}\makebox(0,0)[lb]{\smash{${\cal M}^-$}}}%
    \put(0.76719119,0.21360057){\color[rgb]{0,0,0}\makebox(0,0)[lb]{\smash{${\cal M}^+$}}}%
    \put(0.43066065,0.2947907){\color[rgb]{0,0,0}\makebox(0,0)[lb]{\smash{$u(\cdot,y)$}}}%
    \put(0.87720632,0.41304588){\color[rgb]{0,0,0}\makebox(0,0)[lb]{\smash{$\Gamma$}}}%
  \end{picture}%
\endgroup
\end{center}
Note that the {\it energy} is conserved, indeed if $u\in{\cal H}$ solves $(\ref{eq:eq0})$ on $\R\times (y_{1},y_{2})$ then 
$$
E_{u}(y)= \frac 12\|\partial_{y}u(\cdot,y)\|^{2}_{L^{2}(\R)^{2}}-V(u(\cdot,y))$$
is constant on $(y_{1},y_{2})$ (see \cite{[GuiHamiltonian]} for more general identities of this kind). In particular, for the  heteroclinic type solution $u\in\HH_{m}$ given in Theorem \ref{T:mainPRIMO} we have that $E_{u}(y)=-m$ for every $y\in\R$ and it connects as $y\to\pm\infty$ the two component ${\cal M}_{\pm}$ of the level set $\{q\in\Gamma\,|\,V(q)\leq m\}$.\medskip

\noindent {Now, if we take $c\in(m,m+\lambda)$ with $\lambda>0$ small enough, by $(*)$ we get that
$$
\{q\in\Gamma\,|\, V(q)\le c \}=\G_c^{-}\cup \G_c^{+}\hbox{ with }
\hbox{dist}_{L^{2}}(\G_c^{-},
\G_c^{+})>0.\leqno{(*_{c})}
$$
A natural problem, which generalizes the above one, is to look for  a solution $u\in{\cal H}$ to (\ref{eq:eq0}) with Energy $E_{u}=-c $ which connects in $\Gamma$ the sets $\G_c^{-}$ and $\G_c^{+}$.\par\noindent
In such a case $V(u(\cdot,y))=-E_{u}(y)+\frac 12\|\partial_{y} u(\cdot,y)\|_{L^{2}}^{2}\geq c $ for every $y\in\R$}
and so solutions with energy $-c $ can be sought as minima of the new renormalized functional
\[\varphi_{c}(u)=\int_{\R}\frac 12\|\partial_{y}u(\cdot,y)\|^{2}_{L^{2}(\R)^{2}}+\left(V(u(\cdot,y))-c \right)\,dy\]
on the space
\begin{align*}
{\cal H}_{c}=\{ u\in {\cal H}\,|\, \liminf_{y\to\pm\infty}\hbox{dist}_{L^{2}}(u(\cdot,y),{\cal V}_c^{\pm})=0
\hbox{ and } V(u(\cdot,y))\geq c \hbox{ a.e. in } \R\}.\end{align*}
The functional $\varphi_{c}$ enjoys most of the properties of the functional $\varphi$ and the above concentration arguments  work even in this setting. In particular  the suitable $y$-translated minimizing sequences of $\varphi_{c}$ on ${\cal H}_{c}$ weakly converge to functions $u_{c}\in{\cal H}$ such that $\hbox{dist}_{L^{2}}(u(\cdot,y),{\cal V}_c^{\pm})\to 0$ as $y\to\pm\infty$. However we do not know if it satisfies the constraint $V(u_{c}(\cdot,y))\geq c $ for almost every $y\in\R$ and so that $u_{c}\in\HH_{c}$.
Then, setting
\begin{align*}s_{c}=&\sup\{y\in\R/\hbox{dist}_{L^{2}}(u_{c}(\cdot,y), {\cal V}^{-}_c)\leq d_{0}\hbox{ and }V(u_{c}(\cdot,y))\leq c \}\\
t_{c}=&\inf\{y>s_{c}\, /\,  V(u_{c}(\cdot,y))\leq c \}\end{align*}
what we prove is that 
if $[y_{1},y_{2}]\subset (s_{c},t_{c})$ then $\inf_{y\in[y_{1},y_{2}]}V(u_{c}(\cdot,y))>c$
and so that $u_{c}$ is a solution of (\ref{eq:eq0}) on $\R\times (s_{c},t_{c})$.  By regularity we recover that $V(u(\cdot,y))\to c $ whenever $y\to s_{c}^{+}$ or $y\to t_{c}^{-}$ and the minimality property of $u_{c}$ guarantees that $E_{u_{c}}(y)=-c $ for all $y\in(s_{c},t_{c})$.\smallskip

\noindent In particular, if $s_{c},t_{c}\in\R$, then $V(u(\cdot,t_{c}))=V(u(\cdot,s_{c}))=c $ and hence, by the conservation of energy,  we obtain that
$$
\partial_{y}u(\cdot,t_{c})=\partial_{y}u(\cdot,s_{c})\equiv0 
$$ 
which allows us to recover by reflection from $u_{c}$ a {\it brake orbit type} entire solution. On the other hand, if $s_{c}=-\infty$ (resp. $t_{c}=+\infty$), we can prove that the {$\alpha$-limit} (resp.  {$\omega$-limit}) of $u_{c}$ is constituted by {critical points} of $V$ at level $c$. Then, if $s_{c}=-\infty$ and $t_{c}\in\R$ or  $s_{c}\in\R$ and $t_{c}=+\infty$, from $u_{c}$ we can construct, again by reflection, an {\it homoclinic  type} solution. Finally, if $s_{c}=-\infty$ and $t_{c}=+\infty$ we have that $u_{c}$ is an entire solution of {\it heteroclinic type} connecting $\G^{\pm}_{c}$ as $y\to \pm\infty$
\bigskip
\noindent\def\svgwidth{0.3\columnwidth}

\begingroup
  \makeatletter
  \providecommand\color[2][]{%
    \errmessage{(Inkscape) Color is used for the text in Inkscape, but the package 'color.sty' is not loaded}
    \renewcommand\color[2][]{}%
  }
  \providecommand\transparent[1]{%
    \errmessage{(Inkscape) Transparency is used (non-zero) for the text in Inkscape, but the package 'transparent.sty' is not loaded}
    \renewcommand\transparent[1]{}%
  }
  \providecommand\rotatebox[2]{#2}
  \ifx\svgwidth\undefined
    \setlength{\unitlength}{386.63811035pt}
  \else
    \setlength{\unitlength}{\svgwidth}
  \fi
  \global\let\svgwidth\undefined
  \makeatother
  \begin{picture}(1,0.57596782)%
  \put(0,0){\includegraphics[width=\unitlength]{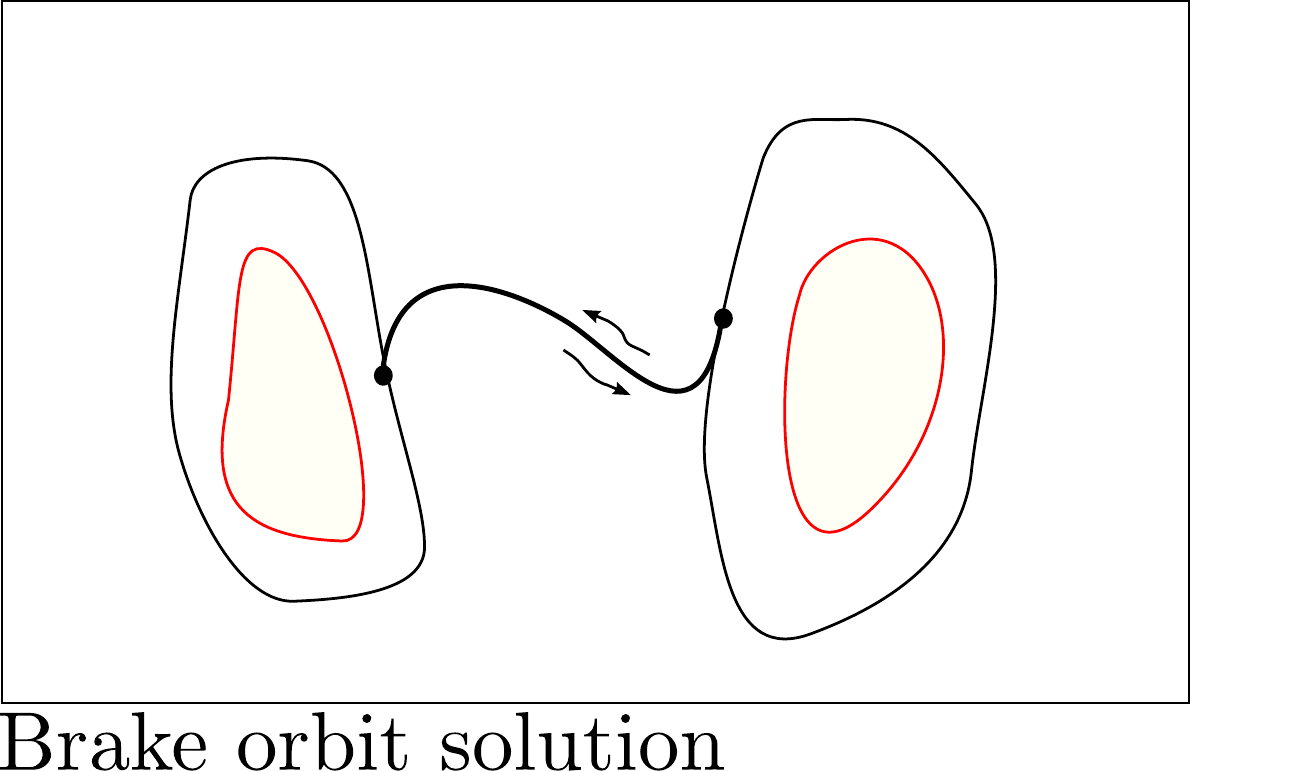}}%
    \put(0.01208382,0.27426216){\color[rgb]{0,0,0}\makebox(0,0)[lb]{\smash{${\cal V}^-_c$}}}%
    \put(0.79588555,0.2846118){\color[rgb]{0,0,0}\makebox(0,0)[lb]{\smash{${\cal V}^+_c$}}}%
    \put(0.82072469,0.50403757){\color[rgb]{0,0,0}\makebox(0,0)[lb]{\smash{$\Gamma$}}}%
    \put(1.1,0){\includegraphics[width=\unitlength]{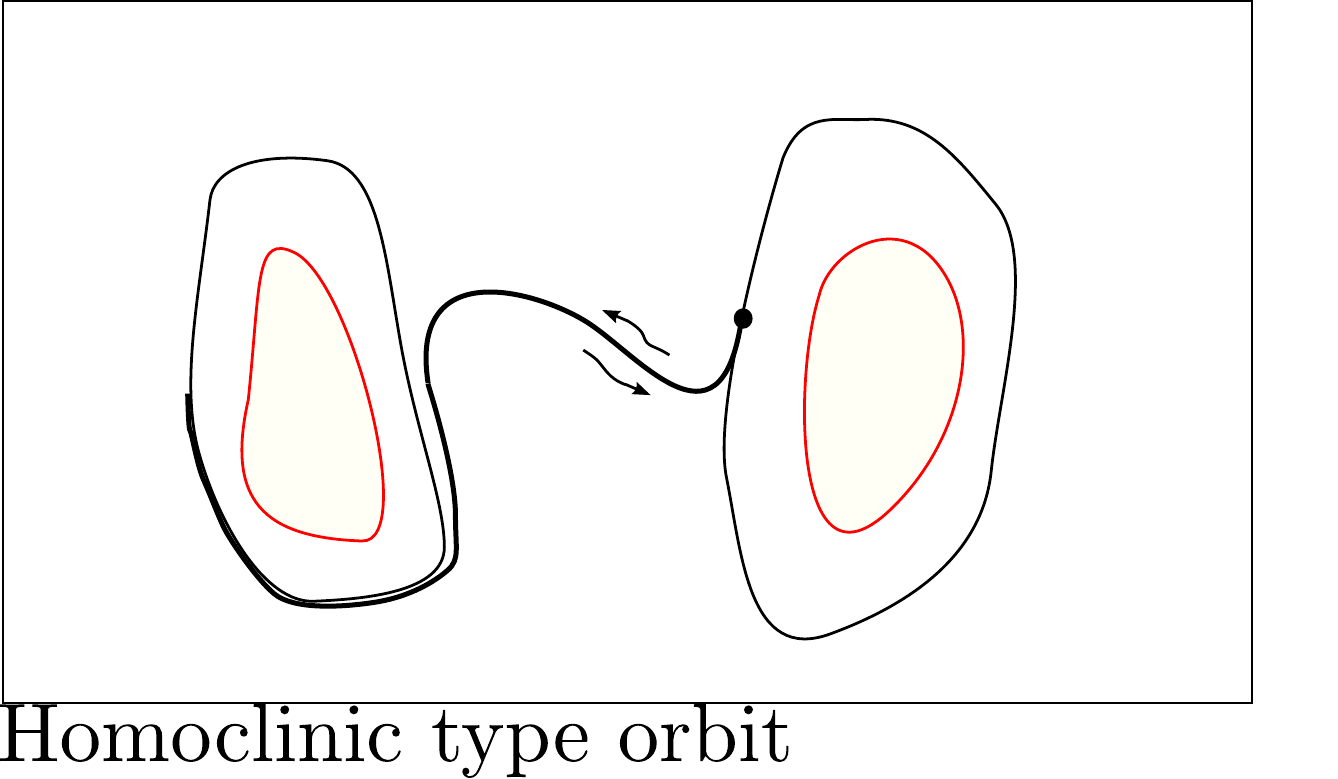}}%
    \put(1.12373167,0.2735746){\color[rgb]{0,0,0}\makebox(0,0)[lb]{\smash{${\cal V}^-_c$}}}%
    \put(1.89428925,0.28069371){\color[rgb]{0,0,0}\makebox(0,0)[lb]{\smash{${\cal V}^+_c$}}}%
    \put(1.96344634,0.50664135){\color[rgb]{0,0,0}\makebox(0,0)[lb]{\smash{$\Gamma$}}}%
    \put(2.2,0){\includegraphics[width=\unitlength]{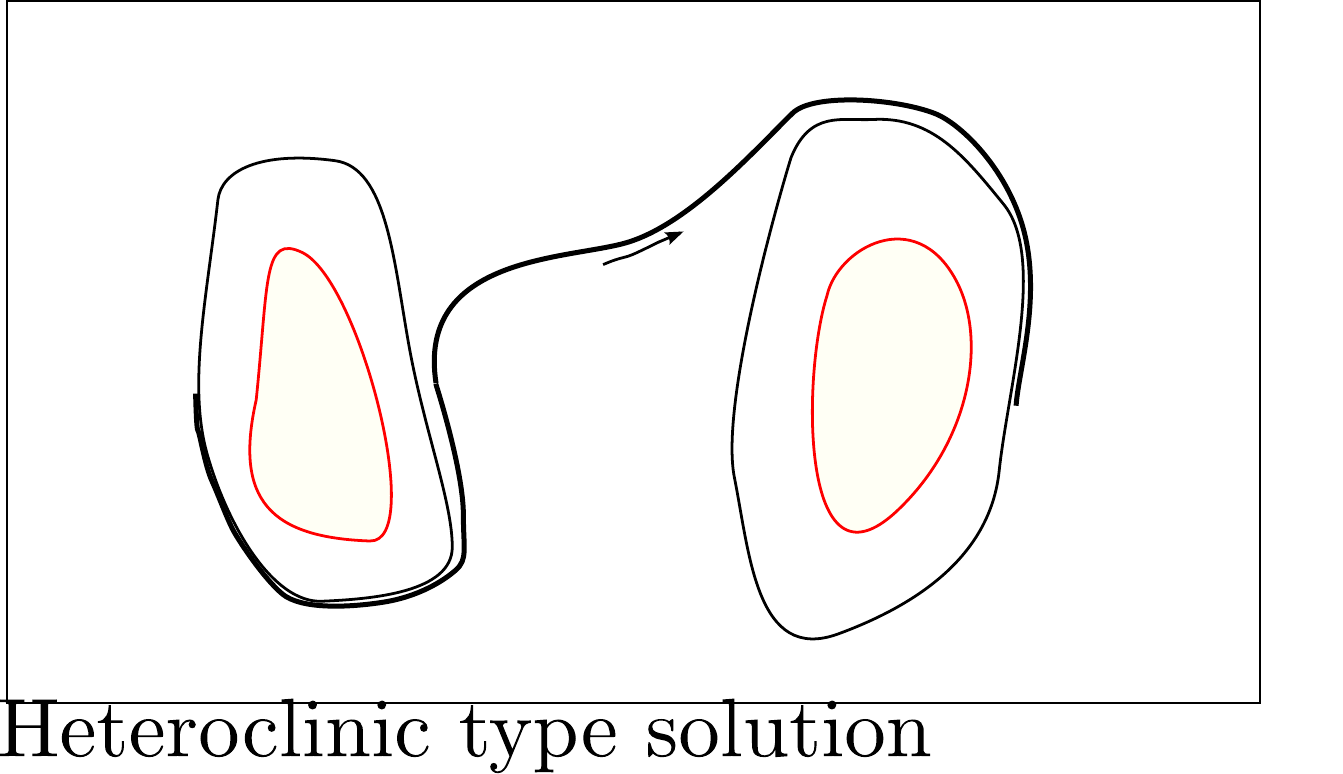}}%
    \put(2.23034729,0.27287713){\color[rgb]{0,0,0}\makebox(0,0)[lb]{\smash{${\cal V}^-_c$}}}%
    \put(2.99850805,0.27388409){\color[rgb]{0,0,0}\makebox(0,0)[lb]{\smash{${\cal V}^+_c$}}}%
    \put(3.06395951,0.50152329){\color[rgb]{0,0,0}\makebox(0,0)[lb]{\smash{$\Gamma$}}}%
  \end{picture}%
\endgroup

\bigskip

\noindent Precisely we prove 

\begin{Theorem}\label{T:main} For every $c \in (m,m+\lambda)$ with $\lambda>0$ small enough, there exists $v_{c}\in\CC^{2}(\R^{2},\R^{2})$ solution to (\ref{eq:eq0})-(\ref{eq:ete}) such that $E_{v_{c}}(y)=-c$ for all $y\in\R$.
Moreover, setting $\KK_{c}^{\pm}=\{q\in\G_c^{\pm}\,|\, V'(q)=0\hbox{ and }V(q)=c \}$, there results
\begin{description}
\item[$(i)$] if $t_{c}=+\infty$ then ${\rm dist}_{H^1(\R)^2}(v_{c}(\cdot,y),\KK_{c}^+)\to 0$ as $y\to+\infty$,
\item[$(ii)$] if $s_{c}=-\infty$ then ${\rm dist}_{H^1(\R)^2}(v_{c}(\cdot,y),\KK_{c}^-)\to 0$ as $y\to-\infty$,
\item[${(iii)}$]  if $t_{c}\in\R$ or $s_{c}\in\R$ then, respectively,  $\partial_{y} v_{c}(\cdot,t_p)\equiv0$ and $v_c(\cdot, t_{c})\in \G_c^+$ or $\partial_{y} v_{c}(\cdot,s_{c})\equiv0$ and $v_{c}(\cdot, s_{c})\in \G_c^-$.
\end{description}
In particular if $c $ is a regular value for $V$ then $t_{c},\, s_{c}\in\R$ and there exists
$T_{c}>0$ such that $v_{c}(x,y+2T_{c})=v_{c}(x,y)$ for all $(x,y)\in\R^{2}$, $\partial_{y} v_{c}(\cdot,0)\equiv\partial_{y} v_{c}(\cdot,T_{c})\equiv0$, $v_{c}(\cdot,0)\in\G_c^{+}$ and $v_{c}(\cdot,T_{c})\in\G_c^{-}$.
\end{Theorem}
\medskip

Note that 
the Theorem guarantees the existence of a brake orbit type solution at level $c$ whenever $c\in(m,m+\lambda)$ is a regular value of $V$. As a consequence of the Sard Smale Theorem and the local compactness properties of $V$, it can be proved that  the set of regular values of $V$ is open and dense in $[m,m+\lambda]$ (see Lemma 2.9 in \cite{[AlMbrake]}). Then, Theorem \ref{T:main} provides in fact the existence of an uncountable set of geometrically distinct two dimensional solutions of (\ref{eq:eq0}) of brake orbit type. If $c\in(m,m+\lambda)$ is a critical value of $V$, then Theorem\ref{T:main} states anyway the existence of an entire solution with energy $-c$ which can be in this case of homoclinic, heteroclinic or brake orbit type, depending on the geometry of the sublevel $\G^{\pm}_{c}$.  \medskip

The variational procedure that we use, and in particular the new renormalized functional $\varphi_c$, was already introduced and used in the framework of scalar non autonomous Allen-Cahn equations in \cite{[AlM3bump]} and \cite{[AlMalmost]}  where we prove the existence of bidimensional solutions that connect as $y\to\pm\infty$ one dimensional solutions which are not minima, in the sense of Giaquinta and Giusti, of the corresponding one dimensional Euler Lagrange functional. Energy prescribed brake orbit  solution were introduced and found in \cite{[AlMbrake]} for the same kind of non autonomous scalar equations.
\smallskip

\section{The one dimensional problem}

In this section  we recall some results concerning the one dimensional heteroclinic problem associated to
(\ref{eq:eq0}),  i.e. the heteroclinic problem 
\begin{equation}\label{eq:ode1}
	\begin{cases}
	-\ddot q(t)+\nabla W(q(t))=0,& t\in\R,\\
	\displaystyle\lim_{t\to\pm\infty}q(t)=\a_{\pm}.& 
	\end{cases}
\end{equation}
Here, we study the problem under the symmetric assumption $(W_{3})$, anyway we want to note that the existence result can in fact be proved without any symmetry of the potential (see e.g. \cite{[Schatzman]} and \cite{[AlikFusco00]}).

Fixed a function $z_{0}\in C^{\infty}(\R,\R^{2})$ 
such that $\|z\|_{\infty}\le R$, $z_{0}$ odd in the first component and even in the second one, $z_{0}(t)=\a_{+}$ for all $t\geq 1$, we consider on the space
$$
\Gamma=z_{0}+H^{1}(\R)^{2},
$$
the functional
$$
V(q)=\int_{\R}\tfrac{1}{2}|\dot q(t)|^{2}+W(q(t))\,dt.
$$
We are interested in the minimal properties of $V$ on $\Gamma$ and we set
$$
m=\inf_{\Gamma}V.
$$ 
Endowing $\Gamma$ with the hilbertian structure induced by the map
$Q:H^{1}(\R)^{2}\to \Gamma$, $Q(z)=z_{0}+z$, it is classical to prove  that $V\in \CC^{2}(\Gamma)$
and that critical points of $V$ are classical solutions to (\ref{eq:ode1}).
\medskip

\noindent Moreover, if $I$ is an interval in $\R$, we set
$$
V_{I}(q)=\int_{I}\tfrac{1}{2}|\dot q(t)|^{2}+W(q(t))\, dt,
$$
noting that $V_{I}(q)$ is  well 
defined on $H^{1}_{loc}(\R)^{2}$ with values in $[0,+\infty]$. \smallskip

\noindent Finally, for a given $q\in L^{2}(\R)^{2}$ we denote $\| q\|\equiv \|
q\|_{L^{2}(\R)^{2}}$ and given $A,\,B\subset L^{2}(\R)^{2}$ we denote
$$
\dd(A,B)=\inf\{\|q_{1}-q_{2}\|\,|\, q_{1}\in A,\, q_{2}\in B\}.
$$

\begin{Remark}\label{R:costanti}
{\rm We precise some basic 
consequences of the assumptions $(W_{1})-(W_{3})$,
fixing some constants. For all $x\in\R^{2}$, we set 
$$
\chi(x)=\min\{|x-\a_{-}|, |x-\a_{+}|\}.
$$
First we note that,  
since $W\in\CC^{2}(\R)$  and $D^{2}W(\a_{\pm})$ are positive definite, then
\begin{equation}\label{BGS}
\hbox{$\forall\, r>0$ $\exists\,\omega_{r}>0$ such that if $\chi(x)\le r$ then }W(x)\geq \omega_{r}\chi(x)^{2}.
\end{equation} 
Then, since $W(\a_{\pm})=0$, $DW(\a_{\pm})=0$ and $D^{2}W(\a_{\pm})$ are positive definite, we have that there exists ${\overline\d}\in (0,\frac{1}{8})$ two constants $\ow>\uw>0$ such that if $ \chi(x)\le 2\overline\d$ then 
\begin{equation}\label{eq:iper}
4\uw|\xi|^{2}\le  D^{2}W(x)\xi\cdot\xi\le4\ow|\xi|^{2} \hbox{ for all }\xi\in\R^{2},
\end{equation}
and
\begin{equation}\label{eq:stimeW}
\uw\chi(x)^{2}\leq W(x)\leq \ow\chi(x)^{2}
 \hbox{ and }|\nabla W(x)|\leq2\uw\chi(x).
\end{equation}

}\end{Remark}

\begin{Remark}\label{R:simm}
{\rm 
Given $q\in\Gamma$, we denote $q(t)=(q(t)_{1},q(t)_{2})$ and we set
$$
t_{-}=\sup\{t\in\R\,/\, q(s)_{1}< 0\,\,\forall s<t\}\quad \hbox{and}\quad
t_{+}=\inf\{t\in\R\,/\, q(s)_{1}> 0\,\,\forall s>t\}.
$$
Now, if $V_{[t_{+},+\infty})(q)\le V_{(-\infty,t_{-}]}(q)$ we set
$$
\hat q(t)=\begin{cases}q(t+t_{+})&\hbox{if $t\geq 0$}\cr
(-q(t_{+}-t)_{1},q(t_{+}-t)_{2})&\hbox{if $t<0$}\end{cases}
$$
while, if $V_{[t_{+},+\infty})(q)> V_{(-\infty,t_{-}]}(q)$, we set
$$
\hat q(t)=\begin{cases}(-q(t_{-}-t)_{1},q(t_{-}-t)_{2})&\hbox{if $t\geq 0$}\cr
q(t+t_{-})&\hbox{if $t<0$}\end{cases}
$$
Then, we obtain that $\hat q\in\Gamma$, $\hat q({0})_{1}=0$ and $\hat q(t)_{1}t>0$ for all $t\not={0}$. Moreover, by definition and the symmetry assumption $(W_{3})$, we obtain that $V(\hat q)\le V(q)$. Indeed, in the first case, by definition and the invariance under translation, we have
$V_{[0,+\infty)}(\hat q)=V_{[t_{+},+\infty)}(q)$
and, by $(W_{3})$, 
\begin{align*}
V_{(-\infty, 0]}(\hat q)&=\int_{0}^{+\infty}\frac12|\dot q(t_{+}-t)|^{2}+W(-q(t_{+}-t)_{1},q(t_{+}-t)_{2})dt\\
&=\int_{t_{+}}^{+\infty}\frac12|\dot q(s)|^{2}+W(-q(s)_{1},q(s)_{2})ds=
V_{[t_{+},+\infty)}(q).
\end{align*}
Hence
$$
V(\hat q)=2V_{[t_{+},+\infty)}(q)\le V_{[t_{+},+\infty)}(q)+V_{(-\infty, t_{-}]}(q)\le V(q)
$$
Analogously in the second case.\\
Moreover note that $\hat q(t)$ is odd in the first component and even in the second one, that is $\hat q\equiv \hat q^{*}$ where we have denoted
$$
p^{*}(t)=(-p(-t)_{1},p(-t)_{2}),\quad\forall t\in\R,\, p\in\Gamma.
$$
In the sequel we will denote 
$$
\Gamma^{*}=\{q\in\Gamma\,|\, q(t)_{1}t>0\,\,\forall t\not=0 \hbox{ and } q\equiv q^{*}\}
$$ 
We have then proved that for all $q\in\Gamma$ there exists $\hat q\in\Gamma^{*}$ such that $V(\hat q)\le V(q)$ and hence that
$$
m=\inf_{\Gamma}V=\inf_{\Gamma^{*}}V.
$$
}\end{Remark}

\begin{Remark}\label{R:taglio}
{\rm We have that for all $\delta\in(0,2\overline\delta)$ if $q\in\Gamma^{*}$ and $|q(t_{0})-\a_{+}|=\delta$ for some $t_{0}>0$ then 
$$
V_{(-\infty,t_{0}]}(q)\geq m-\frac{\delta^{2}}2(1+{2\ow}).
$$
Indeed, let 
$$
\tilde q(t)=\begin{cases}q(t)&\hbox{if $t\le t_{0}$}\cr
(t_{0}+1-t)q(t_{0})+(t-t_{0})\a_{+}&\hbox{if $t_{0}\le t\le t_{0}+1$}\cr
\a_{+}&\hbox{if $t\geq t_{0}+1$}\end{cases}
$$
we have that $\tilde q\in\Gamma$ and therefore
$$
m\le V(\tilde q)=V_{(-\infty, t_{0}]}(q)+\int_{t_{0}}^{t_{0}+1}\frac12|\a_{+}-q(t_{0})|^{2}+W((t_{0}+1-t)q(t_{0})+(t-t_{0})\a_{+})\, dt
$$
But since $\chi((t_{0}+1-t)q(t_{0})+(t-t_{0})\a_{+})\le\delta<2\overline\delta$, by (\ref{eq:stimeW}) we have $W((t_{0}+1-t)q(t_{0})+(t-t_{0})\a_{+})\le\ow\delta^{2}$ and hence we conclude
$$
m\le V_{(-\infty, t_{0}]}(q)+\frac{\delta^{2}}2(1+{2\ow})
$$
Let us fix $\delta_{0}\in(0,\overline\d)$ such that $\frac{\d_{0}}{\overline\d-\d_{0}}<\frac{2\sqrt{2\uw}}{1+2\ow}$ so that 
$$
\lambda_{0}:=\sqrt{2\uw}\d_{0}(\overline\d-\d_{0})-\frac{\d_{0}^{2}}2(1+2\ow)>0.
$$

}\end{Remark}

In the sequel we will study the compactness properties of the sublevel $\{V\le m+\p\}:=\{q\in\Gamma^{*}\,/\,V(q)\leq m+\p\}$.
First of all we note that
if  $q\in H^{1}_{loc}(\R)^{2}$
is such that $W(q(t))\geq \mu$ for all $t\in (\sigma,\tau)\subset\R$, $\mu>0$, then
\begin{equation}\label{eq:stime1dim}
    V_{(\sigma,\tau)}(q)	\geq
	{\textstyle{\frac{1}{ 2(\tau-\sigma)}}}{|q(\tau)-q(\sigma)|}^{2}+
	\mu(\tau-\sigma)
	\geq \sqrt{2 \mu} \
	|q(\tau)-q(\sigma)|.
\end{equation}
As first consequence, using $(W_{2})$ we obtain the following estimate

\begin{Lemma}\label{L:boundLinfty}
      For all $\lambda>0$ there exists $R_{\lambda}>R$ and $C_{\lambda}>0$ such that if $q\in\Gamma^{*}$ and $V(q)\le m+\lambda$  then $\|q \|_{L^{\infty}(\R)^{2}}\leq R_{\lambda}$ and $\|q-z_{0}\|_{H^{1}(\R)^{2}}\le C_{\lambda}$.
\end{Lemma}
\Proof To prove the first assertion, assume by contradiction that there exists $q\in\Gamma^{*}$ such that $V(q)\le m+\lambda$  and $|q(s)|>2r$ for some $r>R$ and $s\in\R$. Then, since $q(t)\to \a_{+}$ as $t\to +\infty$, by continuity we have that there exists 
$s<\sigma<\tau\in\R$ such that $|q(\sigma)|=2r$,
$|q(\tau)|=r$ and $r<|q(t)|<2r$ for any $t\in (\sigma,\tau)$ and hence,
by $(W_{2})$, we obtain that  $W(q(t))\geq\mu_0$ for every $t\in (\sigma,\tau)$. Then, by (\ref{eq:stime1dim}) we conclude
$$
m+\lambda\geq V(q)\geq V_{(\sigma,\tau)}(q)\geq \sqrt{2\mu_0}r
$$
which is impossible for $r>R$ big enough.\\
Now, since $\|q \|_{L^{\infty}(\R)^{2}}\leq R_{\lambda}$, by $(\ref{BGS})$, there exists $\omega_{\lambda}>0$ such that
$$
W(q)\geq \omega_{\lambda}\chi^{2}(q(t)),\quad\forall t\in\R
$$
and since $q(t)_{1}t>0$ for all $t\not=0$, we obtain that $W(q(t))\geq C_{\lambda}|q(t)-\a_{+}|^{2}$ for all $t>0$. Hence
$$
\omega_{\lambda}\int_{0}^{+\infty}|q(t)-\a_{+}|^{2}dt\le \int_{0}^{+\infty}W(q(t)) dt\le V(q)
$$
Therefore, by symmetry, we obtain that $\|q-z_{0}\|_{H^{1}}\le C_{\lambda}$ for some $C_{\lambda}>0$.
\QED

Moreover we have

\begin{Lemma}\label{L:conc}
For all $\lambda\in(0,\lambda_{0})$ there exists $T_{\lambda}>0$ such that if $q\in\Gamma^{*}$ and $V(q)\le m+\lambda$ then $\chi(q(t))< \overline\delta$ for all $|t|>T_{\lambda}$.
\end{Lemma}
\Proof Indeed by  Lemma \ref{L:boundLinfty}, $\|q\|_{L^{\infty}}\le R_{\lambda}$ and by (\ref{BGS}) let $\omega_{\lambda}>0$ be such that $W(q(t))\geq \omega_{\lambda}\chi(q(t))^{2}$ for all $t\in\R$. Then, noting that $q(0)_{1}=0$, $q(t)_{1}>0$ for all $t>0$ and $q(t)\to \a_{+}$ as $t\to +\infty$, 
let $t_{0}=\min\{t>0\,|\, |q(t)-\a_{+}|=\d_{0}\}$. Then, $|q(t)-\a_{+}|\geq\d_{0}$ for all $t\in(0,t_{0})$ and hence for all $t\in(0,t_{0})$ we have
$W(q(t))\geq \omega_{\lambda}\d_{0}^{2}>0$ and then, by (\ref{eq:stime1dim}),  
$$
V(q)\geq V_{(0,t_{0})}(q)\geq \omega_{\lambda}\d_{0}^{2} t_{0}.
$$
Therefore $t_{0}\le T_{\lambda}:=\frac{m+\lambda}{\omega_{\lambda}\delta_{0}^{2}}$.
Now, we claim that $|q(t)-\a_{+}|\le\overline\delta$ for all $t>T_{\lambda}$. Indeed, arguing by contradiction, we have that there exist $\tau>\sigma>T_{\lambda}$ such that $\d_{0}<|q(t)-\a_{+}|<\overline\delta$ for all $t\in(\sigma,\tau)$, $|q(\sigma)-\a_{+}|=\d_{0}$ and $|q(\tau)-\a_{+}|={\overline\delta}$. Hence, by (\ref{eq:stimeW}) we have that $W(q(t))\geq \uw\d_{0}^{2}$ for all $t\in(\sigma,\tau)$ and hence, by (\ref{eq:stime1dim}), 
$$
V_{(\tau,\sigma)}(q)\geq \sqrt{2\uw\delta_{0}^{2}}|q(\tau)-q(\sigma)|\geq\sqrt{2\uw}\d_{0}(\overline\d-\d_{0})
$$ 
But since $|q(t_{0})-\a_{+}|=\d_{0}$ by Remark \ref{R:taglio} we have that
$$
V_{(-\infty,t_{0})}(q)\geq m-\frac{\d_{0}^{2}}2(1+2\ow)
$$
and hence we obtain the contradiction 
$$
m+\lambda\geq V(q)\geq V_{(-\infty,t_{0})}(q)+V_{(\tau,\sigma)}(q)\geq m-\frac{\d_{0}^{2}}2(1+2\ow)+\sqrt{2\uw}\d_{0}(\overline\d-\d_{0})=m+\lambda_{0}.
$$
Then, $|q(t)-\a_{+}|\le \overline\delta$ for all $t>T_{\lambda}$ and, by symmetry, we conclude that $\chi(q(t))\le \overline\delta$ for all $|t|>T_{\lambda}$.
\QED

By Lemma \ref{L:boundLinfty} in particular we  obtain 
\begin{Lemma}\label{L:min}
    Let $(q_{n})$ be a sequence in $\Gamma^{*}$ such that $V(q_{n})\le m+\lambda$ for some $\lambda\in(0,\lambda_{0})$ and all $n\in\N$.
    Then, there exists $q\in
    \Gamma^{*}$ 
    such that, along a subsequence,
    $q_{n}-q\to 0$ weakly in $H^{1}(\R)^{2}$ and $V(q)\le\displaystyle\liminf_{n\to +\infty} V(q_{n})$. 
\end{Lemma}
\Proof By Lemma \ref{L:boundLinfty} we have that there exists $R_{0}>0$ and $C_{0}>0$ such that
$\| q_{n}\|_{L^\infty}\leq R_{0}$ and $\|q_{n}-z_{0}\|_{H^{1}}\le C_{0}$ for all $n\in\N$. Then, up to a subsequence, $(q_{n})$ weakly converge in $\Gamma$ and uniformly on compact subset to some $q\in\Gamma$. Since $(q_{n})\subset\Gamma^{*}$, by uniform convergence we obtain that $q\in\Gamma^{*}$. Finally, by Fatou's Lemma we also obtain that $V(q)\le\displaystyle\liminf_{n\to +\infty} V(q_{n})$.  
  \QED
  
  In particular, using Remark \ref{R:simm}, we obtain that
$\MM=\{q\in\Gamma^{*}\,|\, V(q)=m\}$ is not empty. Then, by standard argument, it can be proved that every $q\in\MM$ is a classical solution to problem (\ref{eq:ode1}).
\bigskip

The following Lemma establishes in particular that $\MM$ is compact.

\begin{Lemma}\label{L:compK}
Let $(q_{n})\subset \Gamma^{*}$ and $q\in\Gamma^{*}$ be such that $q_{n}-q\to 0$ weakly in $H^{1}(\R)^{2}$, $V(q_{n})\to \ell<m+\lambda_{0}$ as $n\to +\infty$ and $V(s_{n}q_{n}+(1-s_{n})q)=\ell$ for all $n\in\N$ and some sequence $(s_{n})\subset[0,1]$. Then, 
    $$
    (1-s_{n})\|q_{n}-q\|^2_{H^{1}}\to 0\hbox{ as }n\to\infty.
    $$
\end{Lemma}
\Proof Since $V(q_{n})\to \ell<m+\lambda_{0}$, we can assume that $V(q_{n})<c+\lambda_{0}$ for all $n\in\N$. Moreover, since $V$ is weakly lower semicontinuous, we have also that $V(q)< m+\lambda_{0}$. Then, by Lemma \ref{L:conc}, there exists $T_{0}>0$ such that 
\begin{equation}\label{eq:conv}
\hbox{if $|t|>T_{0}$ 
then $\chi(q(t))<\overline\d$ and $\chi(q_{n}(t))<\overline\d$ for all $n\in\N$ }
\end{equation}
and the same holds true for their convex combination $s_{n}q_{n}+(1-s_{n})q$. Then, by $(W_{1})$, 
     we obtain that for all $|t|> T_{0}$ there results
     \begin{align} \label{eq:convW1}
     (1-s_{n})W(q_{n}(t))
      \le (1-s_{n})&W(q(t))\\&\nonumber+W(q_{n}(t))-W(s_{n}q_{n}(t)+(1-s_{n})q(t))
     \end{align}
     Moreover, by (\ref{eq:stimeW}),  we have that 
           \begin{equation}\label{eq:convW2}
      {\uw}\,\chi(q_{n}(t))^{2}\le W(q_{n}(t)) \quad\forall |t|> T_{0}.
      \end{equation}
        Since,  by assumption, we have $q_{n}\to q$ in
$L^\infty_{loc}(\R)^{2}$, $\dot q_{n}- \dot q\to0$ weakly  in
$L^2(\R)^{2}$ and 
$V(q_{n})-V(s_{n}q_{n}+(1-s_{n})q)\to 0$,
 we derive that for every $T\geq T_{0}$
\begin{equation}\label{eq:label}
    \frac{1-s_{n}^{2}}{2}\|\dot q_{n}-\dot q\|^{2}+\int_{|t|>T}W(q_{n})-
    W(s_{n}q_{n}+(1-s_{n})q)\, dt\to 0\hbox{ as }n\to\infty.
\end{equation}
Then, by (\ref{eq:convW1}) and (\ref{eq:label}), 
using (\ref{eq:convW2}) we recover 
$$
    \frac{1-s_{n}^{2}}{2}\|\dot q_{n}-\dot q\|^{2}+(1-s_{n})
    {\uw}\int_{|t|>T}\chi(q_{n})^{2}\, dt\leq 
    (1-s_{n})\int_{|t|> T}W(q)\, dt+o(1)\hbox{ as }n\to\infty.
$$ 
Moreover, by the choice of $T_{0}$ and the symmetric property of $q_{n}$ and $q$, for $|t|> T_{0}$ we have
$$
|q_{n}(t)-q(t)|^{2}\leq 2(\chi(q_{n}(t))^{2}+\chi(q(t))^{2}),\quad\forall n\in\N.
$$ 
Hence, since $\int_{\R}W(q)\, dt<+\infty$ and $\int_{\R}\chi(q)^{2}dt<+\infty$, we obtain that for any $\eta>0$ there 
exists $T_{\eta}>T_{0}$ such that 
\begin{eqnarray*}
 && (1-s_{n})(\frac{\|\dot q_{n}-\dot q\|^{2}}{2}+
    {\uw}\int_{|t|>T_{\eta}}|q_{n}- q|^{2}dt)\\
    &&\qquad\le  \frac{1-s_{n}^{2}}{2}\|\dot q_{n}-\dot q\|^{2}+(1-s_{n}){\uw}\int_{|t|>T_{\eta}}|q_{n}- q|^{2}\, 
    dt\leq 
    \eta+o(1)
\end{eqnarray*}
as $n\to +\infty$. Since $q_{n}- q\to 0$ in $L^\infty_{loc}(\R)^{2}$, we obtain the claim.\QED

By Lemma \ref{L:min}, choosing $s_{n}=0$ in Lemma \ref{L:compK}  we can conclude

  \begin{Theorem}\label{T:compK}
  Let $(q_{n})$ be a minimizing sequence for $V$ over $\Gamma^{*}$.
    Then, there exists $q\in
    \MM$ such that, along a subsequence, $\|q_{n}-q\|_{H^{1}}\to 0$ as $n\to +\infty$.
    \end{Theorem}

\noindent In particular, by the previous Lemma we obtain that for every $r>0$ there exists
$\nu_{r}>0$ such that
\begin{equation}\label{eq:h0}
    \hbox{if
    }q\in\{V\leq m+\tfrac\p2\}\hbox{ and }
    \dist_{H^{1}}(q,\MM)\geq r\hbox{
    then }V(q)\geq m+\nu_{r}.
\end{equation}

Now, by condition $(*)$, let $\MM^\pm\subset\MM$ be such that $\MM=\MM^+\cup\MM^-$ and let us denote
\begin{equation}\label{eq:d0}
\dist(\MM^{+},\MM^{-}):=5d_{0}>0.
\end{equation}
Moreover, using (\ref{eq:h0}), we fix $m^{*}\in(m,m+\frac{\lambda_{0}}2)$ such that
\begin{equation}\label{eq:c*}
\hbox{ if }q\in\{V\leq m^{*}\}\hbox{ then }\dist_{H^{1}}(q,\MM)\leq d_{0},
\end{equation}
then, we set
$$
\G^{\pm}=\{q\in\Gamma^{*}\,|\,V(q)\le m^{*}\hbox{ and } \dist_{H^{1}}(q,\MM^{\pm})\leq d_{0}\}.
$$
Plainly we have that $\MM^{\pm}\subset\G^{\pm}$ and, by (\ref{eq:d0}) and (\ref{eq:c*}), we have
\begin{equation}\label{eq:dis}
\{V\le m^{*}\}=\G^{+}\cup\G^{-}\hbox{ and }\dist(\G^{+},\G^{-})\geq 3d_{0}
\end{equation}
Finally, fixed any $c\in[m,m^*]$, let us denote $\G_c^\pm=\{u\in\G^\pm\,|\, V(u)\le c\}$ and note that 
$$
\hbox{$\{V\le c\}=\G^+_c\cup\G^-_c$ and $\dist(\G_c^{+},\G_c^{-})\geq 3d_{0}$.}
$$
Moreover, we have

\begin{Lemma}\label{R:chiusuraG+}
The sets $\G_c^{\pm}$ are weakly compact in $\Gamma^*$.
\end{Lemma}
{\Proof Indeed, if $(q_{n})\subset\G_c^{\pm}$ then $\dist_{H^{1}}(q_{n},\MM^{\pm})\leq d_{0}$ and $V(q_{n})\leq m^{*}$.
Then, by Lemma \ref{L:min}, there exists $q\in\{V\le c\}$ such that, up to a subsequence, $q_{n}-q\to 0$ weakly in $H^{1}(\R)^{2}$. Hence, by the weak semicontinuity of the norm, we obtain $\inf_{\bar q\in\MM^{\pm}}\|q-\bar q\|_{H^{1}}\leq d_{0}$ and so $q\in\G^{\pm}$.}\QED

\section{Two dimensional solutions}

In this section we will prove Theorems \ref{T:mainPRIMO} and \ref{T:main} using a variational approach. First,
for $(y_{1},y_{2})\subset\R$ we set
$S_{(y_{1},y_{2})}=\R\times (y_{1},y_{2})$ and
we consider the space
$$
\SS=\{ u\in H^{1}_{loc}(\R^{2})^{2}\,|\,u-z_{0}\in
\cap_{(y_{1},y_{2})\subset\R}H^{1}(S_{(y_{1},y_{2})})^{2}\}.
$$
Note that if $u\in\SS$ then,  $u(\cdot,y)\in\Gamma$ for almost every 
$y\in\R$ and we will consider
$$
\SS^{*}=\{u\in\SS\,|\, u(\cdot,y)\in\Gamma^{*}\hbox{ for almost every }
y\in\R\}.
$$
Moreover note that, setting
 $$\overline{\Gamma}=z_{0}+L^{2}(\R)^{2},$$ 
 the completion
of $\Gamma$ with respect to the $L^{2}$-metric, for all $u\in\HH$ the function $y\in\R\mapsto u(\cdot,y)\in \overline\Gamma$ 
defines a continuous trajectory verifying
\begin{equation}\label{eq:continuita1}
    \|u(\cdot, y_{2})-u(\cdot, y_{1})\|^{2} \leq
    \|\partial_{y} u\|^{2}_{L^{2}
    (S_{(y_{1},y_{2})})^{2}}
    |y_{2}-y_{1}|,\quad\forall\, (y_{1},y_{2})\subset\R.
\end{equation}

\begin{Remark}\label{R:continuita}
{\rm  Given any $u\in\SS^{*}$
and $(y_{n})\subset\R $ such that $u(\cdot, y_{n})\in\G_c^{\pm}$ and $y_{n}\to y_{0}$ as $n\to +\infty$ then, by (\ref{eq:continuita1}), $u(\cdot,y_{n})\to u(\cdot,y_{0})$ in $\overline\Gamma$ and since, by Lemma\ref{R:chiusuraG+}, $\G_c^{\pm}$ are weakly precompact in $\Gamma^*$, we conclude that $u(\cdot,y_{n})\to u(\cdot,y_{0})$ weakly in $\Gamma^{*}$ and $u(\cdot, y_{0})\in\G_c^{\pm}$.
}\end{Remark}

Finally, setting $V(u(\cdot,y))=+\infty$ whenever $u(\cdot,y)\in\overline\Gamma\setminus\Gamma$, we note that the function $y\in\R\mapsto V(u(\cdot,y))\in [m,+\infty]$ is lower semicontinuous for every given $u\in\HH$.\bigskip

\noindent
Fixed  a value $c \in [m,m^{*})$,  we look for minimal 
properties of the functional
$$
\f(u)=\int_{\R}
\tfrac{1}{2}\|\partial_{y}u(\cdot,y)\|^2+
	(V(u(\cdot,y))-c ) \,dy.
$$
on the set
\begin{align*}
{\cal H}_{c}=\{ u\in {\cal H}^*\,|\, \liminf_{y\to\pm\infty}\dist(u(\cdot,y),{\cal V}_c^{\pm})\le d_0
\hbox{ and } V(u(\cdot,y))\geq c \hbox{ a.e. in } \R\}.\end{align*}
Note that if $u\in\HH_{c}$ then 
    $V(u(\cdot,y))\geq c $
    for almost every $y\in\R$ and so the functional $\f$ is well defined on $\HH_{c}$
    with values in $[0,+\infty]$.\\
 Finally, given an interval $I\subset\R$ we also consider on $\HH_{c}$ the functional
    $$
	\varphi_{c,I}(u)=\int_{I}\tfrac{1}{2}\|\partial_{y}u(\cdot,y)\|+
	    (V(u(\cdot,y))-c ) \,dy
    $$

\begin{Remark}\label{R:Mnonvuoto}{\rm
    \begin{itemize}
        \item[$(i)$] Given any $I\subseteq \R$,  the functional $\varphi_{c,I}(u)$  is well defined for all $u\in \SS^{*}$
    such that the set $\{ y\in I\, |\, V(u(\cdot,y))<c \}$
    has bounded measure. Moreover if $(u_{n})\subset\HH_{c}$ is such that $u_{n}\to u$ weakly in
    $H^{1}_{loc}(\R^{2})^{2}$, then
    $\varphi_{c,I}(u)\leq{\liminf}
    \varphi_{c,I}(u_{n})$ (see \cite{[AJM]}, Lemma 3.1).
    
        \item[$(ii)$]  \ Given $u\in \SS^{*}$, if $(y_{1},y_{2})\subset\R$ and $\mu>0$ are
    such that
    $V(u(\cdot,y))\geq c +\mu$
    for all $y\in (y_{1},y_{2})$, then
    \begin{align}	 
	\varphi_{c,(y_{1},y_{2})}(u)
		  \label{eq:stime2dim}
		&\geq\tfrac{1}{2(y_{2}-y_{1})}
		\| 
		u(\cdot,y_{1})-u(\cdot,y_{2})\|^{2}+\mu(y_{2}-y_{1})\\
		\nonumber&\geq 
		\sqrt{2\mu}\,\| u(\cdot,y_{1})-u(\cdot,y_{2})\|.
	\end{align}
	 \item[$(iii)$] There result $\HH_{c}\not=\emptyset$ and, by (\ref{eq:stime2dim}), $\inf_{\HH_{c}}\f:= m_{c}\geq\sqrt{2(m^{*} -c )}d_{0}$.
	\end{itemize}
	    }
\end{Remark}

\noindent We 
characterize here below some 
concentration and compactness properties of the minimizing sequences in $\HH_{c}$. First, thanks to the discreteness of the set $\{ V\leq m^{*}\}$ described by (\ref{eq:dis}) we obtain
that if
$\f(u)<+\infty$ 
then, the trajectory
$y\in\R\to u(\cdot,y)\in \overline\Gamma$ is bounded.
Precisely

\begin{Lemma}\label{L:sceltap} There exists $\MC>0$ such that
 if $u\in\HH_{c}$  then
$\| u(\cdot,y_{1})-u(\cdot,y_{2})\|_{L^{2
}}\leq
\MC\, \f(u)$ for all $(y_{1},y_{2})\subset\R$.
\end{Lemma}
\Proof First note that, by Lemma \ref{R:chiusuraG+},
$\diam(\G^{+})$ and $\diam(\G^-)$ are bounded, let us denote $\MD$ the maximum value. Then, the Lemma simple follows noting that $\{V\le m^*\}=\G^+\cup\G^-$ and  if $V(u(\cdot,y))>m^{*}>c$ for all $y\in(\sigma,\tau)\subseteq(y_{1},y_{2})$, by (\ref{eq:stime2dim})
we have
$$
\| u(\cdot,\sigma)-u(\cdot,\tau)\|\leq
\frac{\varphi_{p}(u)}{\sqrt{2(m^{*}-c)}}.
$$
Hence, for all $y\in[y_1,y_2]$ we have that $\dist(u(\cdot,y),\G^+\cup\G^-)\le \frac{\varphi_{p}(u)}{\sqrt{2(m^{*}-c)}}$ and then
$$
\| u(\cdot,y_{1})-u(\cdot,y_{2})\|\leq
2\MD+
\frac{\varphi_{p}(u)}{\sqrt{2(m^{*}-c)}}
$$
and the lemma follows choosing $\MC=\max\{\frac1{\sqrt{2(m^{*}-c )}};2\frac{\MD}{m_{c}}\}$.\QED

By the previous result we get the following characterization of the asymptotic behaviour of the trajectory 

\begin{Lemma}\label{L:liminfinito}
    If  $u\in\HH_{c}$ and
    $\f(u)\leq m_{c}+1$ then
    $
     \displaystyle\lim_{y\to\pm\infty}\dist(u(\cdot,y),\G_c^{\pm})=0$.\end{Lemma}
\Proof Just note that
if $u\in\HH_{c}$ and $\f(u)\leq m_{c}+1$
then $\liminf_{y\to\pm\infty} V(u(\cdot,y))=c $ and so
$\liminf_{y\to\pm\infty}\dist(u(\cdot,y),\{V\le c\})=0$. Since $\{V\le c\}=\G^{-}_c\cup\G_c^{+}$ and $\dist(\G_c^{-},\G_c^{+})\geq 3d_{0}$, by the definition of $\HH_{c}$ we get 
$$
\displaystyle\liminf_{y\to\pm\infty}\dist(u(\cdot,y),\G_c^{\pm})=0.
$$
Now, by contradiction, assume that $
\limsup_{y\to+\infty}\dist(u(\cdot,y),\G_c^+)>0$. Then, by (\ref{eq:continuita1}) the path $y\mapsto u(\cdot,y)$
crosses infinitely many times an annulus of positive thickness $d>0$ around $\G^+_c$ in the $L^{2}$ metric. Then, by
(\ref{eq:stime2dim}), we get the contradiction $\varphi_c(u)=+\infty$. Analogously we obtain $
\limsup_{y\to-\infty}\dist(u(\cdot,y),\G_c^-)=0$\QED

To recover compactness properties, modulo $y$-translations, of the minimizing sequences of $\f$ on
$\HH_{c}$, it is useful to study the concentration properties of the trajectories $y\in\R\to
u(\cdot,y)\in\overline\Gamma$ when $\f(u)$ is close to $m_{c}$. 
To this aim, by Lemmas \ref{L:min} e \ref{L:compK},  we have first the following result

\begin{Lemma}\label{L:tagliofond}
    Let $(u_{n})\subset\HH_{c}$ and $y_{n}\in\R$ be such that  
    $V(u_{n}(\cdot,y_{n}))\to c $ as $n\to +\infty$. Then
    \begin{itemize}
    \item[$(i)$] if  $u_{n}(\cdot,y_{n})\in\G^{+}$ for  all $n\in\N$  then, 
    $ \displaystyle\liminf_{n\to+\infty}\varphi_{c,(-\infty,y_{n})}(u_{n})\geq m_{c}$
    \item[$(ii)$] if 
    $u_{n}(\cdot,y_{n})\in\G^{-}$ for  all $n\in\N$ then, 
    $\displaystyle\liminf_{n\to+\infty}\varphi_{c,(y_{n},+\infty)}(u_{n})\geq m_{c}$.
\end{itemize}
\end{Lemma}
\Proof  { We give only a sketch of the proof, the details can be found in the proof Lemma 3.4 in \cite{[AlMalmost]}}.
Let us prove $(i)$. Note that, by the invariance with respect to the 
    $y$-translation of $V$, we can assume $y_{n}=0$ for all $n\in\N$. 
Setting 
$q_{n}=u_{n}(\cdot,0)$, there results $q_{n}\in\G^{+}$ and $V(q_n)\to c$, hence,
by Lemmas \ref{L:min} and \ref{R:chiusuraG+}, there exists 
$q\in \G^{+}$
such that $V(q)\leq c $ and, up to a subsequence, $q_{n}-q\to 0$ 
weakly in $H^{1}(\R)^{2}$.
We set $s_{n}=\sup\{s\in[0,1]\,|\,
V(sq_{n}+(1-s)q)\leq c \}$ nothing that,
 by continuity,
$V(s_{n}q_{n}+(1-s_{n})q)=c $ for all $n\in\N$. Hence, by Lemma \ref{L:compK}, we have 
\begin{equation}\label{eq:tnvn1}
    (1-s_{n})\|q_{n}-q\|^{2}_{H^{1}}\to 0\hbox{ as }n\to\infty.
\end{equation}
Now,  considering
the new sequence
$$\tilde u_{n}(x,y)=\begin{cases}u_{n}(x,y)&\hbox{if }y\leq 0,\\
			   yq(x)+(1-y)q_{n}(x)& 
			   \hbox{if }0\leq
			   y\leq 1-s_{n},\\
			  (1-s_{n})q(x)+ s_{n}q_{n}(x)& \hbox{if }y\geq 1-s_{n},
			   \end{cases}
$$
we have $\tilde u_{n}\in\HH_{c}$ and so 
$\f(\tilde u_{n})\geq m_{c}$. Then, since
$$
\varphi_{c,(-\infty,0)}(u_{n})=\f(\tilde u_{n})
-\int_{0}^{1-s_{n}}
    \tfrac{1}{2}\|q_{n}-q\|_{L^{2}}^{2}+
    (V(yq(x)+(1-y)q_{n}(x))-c) \, dy,$$
the Lemma follows once we prove that
\begin{equation*}\label{eq:claim2}
    \int_{0}^{1-s_{n}}
    \tfrac{1}{2}\|q_{n}-q\|_{L^{2}}^{2}+
    V(yq(x)+(1-y)q_{n}(x))-c \, dy\to 0,\quad\hbox{as $n\to+\infty$.}
\end{equation*}
Indeed,
$\int_{0}^{1-s_{n}}\frac{1}{2}\|q_{n}-q\|^{2}\, dy=
\frac{1}{2}(1-s_{n})\|q_{n}-q\|^{2}_{L^{2}}\to 
0$ follows by (\ref{eq:tnvn1}). 
While, by (\ref{BGS}), as in \cite{[AlM3bump]}, Lemma 2.13, we have that 
there exists $C>0$  such that
\begin{align*}
V(yq+(1-y)q_{n})-c &=V(yq+(1-y)q_{n})-V((1-s_{n})q+s_{n}q_{n})\\
&
\leq C(1-s_{n})\|q_{n}-q\|_{H^{1}}
\end{align*}
for all $n\in\N$ and $y\in [0,1-s_{n}]$. Then, by 
(\ref{eq:tnvn1}), $\int_{0}^{1-s_{n}}
V(yq+(1-y)q_{n})-c \, dy\to 0$ and the Lemma follows.
\QED

Lemma \ref{L:tagliofond} is used to obtain the following result
that, together with (\ref{eq:stime2dim}), is used in Lemma \ref{L:primosol} 
to characterize the concentration properties
of the minimizing sequences of $\f$ on $\HH_{c}$.
First, denoting
$$
\ell_{0}=\min\{1;\sqrt{\tfrac{m^{*}-c }2}d_{0}\},
$$
we obtain the following concentration result 
\begin{Lemma}\label{L:cattura}
    There exists $\nu\in (0,m^{*}-c )$ such that
	if $u\in\HH_{c}$,
	$\f(u)\leq m_{c}+\ell_{0}$
	and $V(u(\cdot,\bar y))< c +\nu$ for some 
	$\bar y\in\R$, then,
	either
	\begin{itemize}
	    \item[$(i)$] $u(\cdot,\bar y)\in
	\G^{-}$  and $\dist(u(\cdot,y),
	\G^{-})\leq d_{0}$ for all $y\leq\bar y$; or

	    \item[$(ii)$] $u(\cdot,\bar y)\in
	    \G^{+}$ and $\dist(u(\cdot,y),
	    \G^{+})\leq d_{0}$ for all $y\geq\bar y$.
	\end{itemize}
\end{Lemma}
\Proof 
We prove only $(ii)$ since the proof of $(i)$ is analogous.\\
Assume
by contradiction that there exist 
a sequence
$(u_{n})\subset\HH_{c}$ such that $\f(u_{n})\leq 
m_{c}+\ell_{0}$, 
and two sequences $y_{n,1}<y_{n,2}$ in $\R$ such that for all
$n\in\N$ there
results
$$  
    \lim_{n\to\infty} V(u_{n}(\cdot,y_{n,1}))= c ,\ 
     u_{n}(\cdot,y_{n,1})\in\G^{+}\,\hbox{ and }\,
     \dist(u_{n}(\cdot,y_{n,2}),\G^{+})>d_{0}.
$$
Since $u_{n}(\cdot,y_{n,1})\in\G^{+}$ and
$\dist(u_{n}(\cdot,y_{n,2}),\G^{+})>d_{0}$, by (\ref{eq:continuita1}), since $\dd(\G^{+},\G^{-})\geq 3d_{0}$ and $\G^-\cup\G^+=\{V\le m^*\}$,
we obtain that there exists $\bar y_{n,1},\bar y_{n,2}\in [y_{n,1},y_{n,2}]$
such that $V(u_{n}(\cdot,y))\geq m^{*}$ for every $y\in (\bar y_{n,1},\bar
y_{n,2})$
and $\dist(u_{n}(\cdot,\bar y_{n,1}),u_{n}(\cdot,\bar y_{n,2}))=d_{0}$. By
(\ref{eq:stime2dim})
we obtain
$\varphi_{c,(y_{n,1},+\infty)}(u_{n})\geq
\varphi_{c,(\bar y_{n,1},\bar y_{n, 2})}(u_{n})\geq
\sqrt{2(m^{*}-c )}d_{0}\geq 2\ell_{0}$
for all $n\in\N$.
Hence, for all $n\in\N$ we get
$$
\varphi_{c,(-\infty,y_{n,1})}(u_{n})=\varphi_{c}(u_{n})-\varphi_{c,(y_{n,1},+\infty)}(u_{n})\le m_{c}-\ell_{0}
$$
which is a contradiction since, by Lemma \ref{L:tagliofond}, $\liminf\varphi_{c,(-\infty,y_{n,1})}(u_{n})\geq m_{c}$.
\QED

Hence, using Lemmas \ref{L:sceltap}, \ref{L:liminfinito} and \ref{L:cattura}, arguing as in Lemma 3.6 in \cite{[AlMbrake]}, we can prove

\begin{Lemma}\label{L:primosol}
    Let $(u_{n})\subset\HH_{c}$
    be such that $\f(u_{n})\to m_{c}$ as $n\to\infty$ and
    such that $\dist(u_{n}(\cdot,0),\G^{-})=d_{0}$ for all
    $n\in\N$.
    Then, there exists $u_{c}\in \HH^{*}$ such that, up to a subsequence,
    \begin{itemize}
	\item[$(i)$] $u_{n}-u_{c}\to 0$ as $n\to\infty$
    weakly in $H^{1}_{loc}(\R^{2})^{2}$,


	\item[$(ii)$] there exists $L_{c}>0$ such that 
$\dist(u_{c}(\cdot,y),\G^{-})\leq d_{0}$ for all $y\leq -L_{c}$,  and 
$\dist(u_{c}(\cdot,y),\G^{+})
\leq d_{0}$ for all $y\geq L_{c}$.
    \end{itemize}
\end{Lemma}
\Proof
Since $\f(u_{n})\to m_{c}$
as $n\to\infty$ we can assume that $\f(u_{n})\leq 
m_{c}+\ell_{0}$ for all
$n\in\N$. To prove $(i)$ we show that there exists $u_{c}\in
H^{1}_{loc}(\R^{2})$  such that, along a subsequence,
$u_{n}-z_{0}\to u_{c}-z_{0}$ weakly 
in $H^{1}(S_{(-k,k)})^{2}$ for every $k\in\N$. This plainly implies that $u_{c}\in\HH^{*}$ and that $u_{n}-u_{c}\to 0$ weakly in $H^{1}_{loc}(\R^{2})^{2}$.\\
To this aim first note that fixed any function $q\in\G^{-}$,  
by Lemmas \ref{L:sceltap} e \ref{L:liminfinito}, since $\dd(\G^{+},\G^{-})\geq 3d_{0}$, we have
that there exists $C>0$ such that for all $y\in\R$ there results
$$
\|u_{n}(\cdot,y)-q\|\leq\dist(u_{n}(\cdot,y),\G^{-})+\diam(\G^{-})\leq
C.
$$
Then,  $\|u_{n}-q\|^{2}_{L^{2}(S_{(-k,k)})^{2}}\leq 
2kC^{2}$ for all $n\in\N$ and $k\in\N$. Since moreover
$\|\nabla u_{n}\|^{2}_{L^{2}(S_{(-k,k)})^{2}}\leq
2(\f(u_{n})+2kc )$ we conclude that the sequence
$(u_{n}-q)_{n\in\N}$, and so
the sequence $(u_{n}-z_{0})_{n\in\N}$,
is
bounded in $H^{1}(S_{(-k,k)})^{2}$ for every $k\in\N$. Then, a diagonal argument
implies the existence of a function $u_{c}\in
H^{1}_{loc}(\R^{2})^{2}$  such that along a subsequence
$u_{n}\to u_{c}$ weakly 
in $H^{1}(S_{(-k,k)})^{2}$ for every $k\in\N$ and $(i)$ follows.\\ 
Finally, to prove $(ii)$, note that by  (\ref{eq:stime2dim}) there exists
$L_{c}>0$ such that, for every $n\in\N$ there exist $y_{n,1}\in (-L_{c},0)$ and
$y_{n,2}\in (0,L_{c})$ such that $V(u_{n}(\cdot,y_{n,1}))$, 
$V(u_{n}(\cdot,y_{n,2}))\leq c +\nu$. By Lemma \ref{L:cattura}, since $\dist(u_{n}(\cdot,0),\G^{-})=d_{0}$, we obtain that for
every $n\in\N$ there results
$$
u_{n}(\cdot,y_{n,1})\in\G^{-}\hbox{
and }\dist(u_{n}(\cdot,y),\G^{-})\leq d_{0}\hbox{ for all }y\leq-L_{c}
$$
and 
$$
u_{n}(\cdot,y_{n,2})\in
    \G^{+}\hbox{ and }\dist(u_{n}(\cdot,y),\G^{+})\leq d_{0}\hbox{ for all }y\geq
L_{c}.
$$
Hence, in the limit we obtain
$$
\dist(u_{c}(\cdot,y),\G^{-})\leq d_{0}\hbox{ for all }y\leq-L_{c}\hbox{
and }\dist(u_{c}(\cdot,y),\G^{+})\leq d_{0}\hbox{ for all }y\geq
L_{c}$$
and $(ii)$ follows.\QED

    By the invariance with respect to the 
    $y$-translation of $\f$ and the definition of $\HH_{c}$, we have that there exists a minimizing sequence  $(u_{n})$ which verifies the condition
    $\dist(u_{n}(\cdot,0),\G^{-})=d_{0}$ for all $n\in\N$. Then, by 
Lemma \ref{L:primosol}, the sequence weakly converge in $H^{1}_{loc}(\R^{2})^{2}$ to a function $u_{c}\in\HH^{*}$ such that
    \begin{equation*}\label{eq:minimo}
   \hbox{$\dist(u_{c}(\cdot,y),\G^{-})\leq d_{0}$ for all $y\leq -L_{c}$,  and 
$\dist(u_{c}(\cdot,y),\G^{+})
\leq d_{0}$ for all $y\geq L_{c}$}
    \end{equation*}
In the sequel we will study the minimality properties of the limit point $u_{c}$ which will be used to recover from it a solution $v_{c}$ to (\ref{eq:eq0}) such that $v_{c}(x,y)\to \a_{\pm}$ as $x\to \pm\infty$ uniformly w.r.t. $y\in\R$.\medskip

First of all we remark that if $c =m$ then $u_{m}\in\HH_{m}$, indeed $u_{m}(\cdot,y)\in\Gamma^{*}$ for almost every $y\in\R$ and the condition $V(u_{m}(\cdot,y))\geq m$ is satisfied for almost every $y\in\R$. On the other hand, if $c >m$ 
we do not know if it satisfies the constraint $V(u_{c}(\cdot,y))\geq c $ for almost every $y\in\R$ and then that $u_{c}\in\HH_{c}$. Anyhow
 we will prove that such condition is satisfied on the interval $(s_{c},t_{c})$ where $s_{c}$ and $t_{c}$ are defined as follows:
$$s_{c}=\sup\{y\in \R\, /\, \dist(u(\cdot,y),\G_c^{-})\leq d_{0}\hbox{ and }
			V(u(\cdot,y))\leq
			c \},
$$
$$t_{c}=\inf\{y>s_{c}\, /\, 
			V(u(\cdot,y))\leq
			c \},
$$
where we agree that $s_{c}=-\infty$ whenever $V(u_{c}(\cdot,y))>c $
for every $y\in\R$ such that $\dist(u_{c}(\cdot,y),\G_c^{-})\leq d_{0}$ and that
$t_{c}=+\infty$ whenever $V(u_{c}(\cdot,y))>c $
for all $y>s_{c}$. 

\begin{Remark}\label{R:tpm}{\rm Note that if $c =m$ then $s_{c}=-\infty$ and $t_{c}=+\infty$. While if $s_{c}, t_{c}\in\R$ we have that $s_{c}< t_{c}$, $s_{c}\leq L_{c}$ and $t_{c}\geq -L_{c}$. 
}
\end{Remark}

\medskip

We can now display the minimality properties of the function $u_{c}$.

\begin{Lemma}\label{L:taglio2(i)} 
For every $[y_{1},y_{2}]\subset (s_{c},t_{c})$ there results $\inf_{y\in[y_{1},y_{2}]}V(u_{c}(\cdot,y))>c $.  Moreover,  $\varphi_{c,(s_{c},t_{c})}(u_{c})\le m_{c}$.
\end{Lemma}
\Proof Note that by definition of  $s_{c}$ and $t_{c}$, we have that $V(u_{c}(\cdot,y))>c $ for any $y\in (s_{c},t_{c})$. Then, since the function $y\mapsto V(u_{c}(\cdot,y))$ is semicontinuous, we derive that $\inf_{y\in[y_{1},y_{2}]}V(u_{c}(\cdot,y))>c $ whenever $[y_{1},y_{2}]\subset (s_{c},t_{c})$. 
Finally, 
 $\varphi_{c,(s_{c},t_{c})}(u_{c})$ is well defined and $\varphi_{c,(s_{c},t_{c})}(u_{c})\leq m_{c}$ follows  by
 definition of $u_{c}$ and Remark \ref{R:Mnonvuoto}-$(i)$.\QED
 
 In particular, by Lemmas \ref{L:taglio2(i)}  and \ref{L:primosol}-$(ii)$, arguing as in Lemma \ref{L:liminfinito}, if $s_{c}=-\infty$ then 
$$
\displaystyle\lim_{y\to -\infty}\dd(u_{c}(\cdot,y),\G_c^{-})=0\hbox{ and }\displaystyle\liminf_{y\to -\infty}V(u_{c}(\cdot,y))=c .
$$ 
 Analogously, if $t_{c}=+\infty$ then 
 $$
 \displaystyle\lim_{y\to +\infty}\dd(u_{c}(\cdot,y),\G_c^{+})=0\hbox{ and }\displaystyle\liminf_{y\to +\infty}V(u_{c}(\cdot,y))=c .
 $$

On the other hand, if $s_{c}\in\R$ or $t_{c}\in\R$ we obtain
 
 \begin{Lemma}\label{L:taglio2(ii)} 
If $s_{c}\in\R$ then
$u_{c}(\cdot,s_{c})\in\G_c^{-}$ 
and analogously,  if $t_{c}\in\R$ then
$u_{c}(\cdot,t_{c})\in\G_c^{+}$. 
\end{Lemma}
\Proof If $s_{c}\in\R$, by definition there exists a sequence  such that $y_{n}\to s_{c}^{-}$ as $n\to+\infty$, $V(u_{c}(\cdot,y_{n}))\leq c $ and $\dist(u_{c}(\cdot,y_{n}),\G_c^{-})\leq d_{0}$ for every $n\in\N$. Then, since $\{V\le c\}=\G_c^{-}\cup\G_c^{+}$  we obtain that $u_{c}(\cdot,y_{n})\in\G_c^{-}$ for every $n\in\N$. Moreover, by Remark \ref{R:continuita}, we get $u_{c}(\cdot,y_{n})\to u_{c}(\cdot,s_{c})$ weakly in $\Gamma$ and $u_{c}(\cdot,s_{c})\in\G_c^{-}$. 
\QED

Moreover

\begin{Lemma}\label{L:taglio2}
{\rm If $s<t\in\R$ and
$u\in\HH^{*}$ verify $u(\cdot,s)\in\G_{c}^{-} $,
$u(\cdot,t)\in\G_c^{+}$ 
and $V(u(\cdot,y))\geq c $ for every $y\in (s,t)$,
  then $\varphi_{c,(s,t)}(u)\geq m_{c}$.}
\end{Lemma}
{\Proof Let us fix  two sequences $(s_{n}),\, (t_{n})\subset (s,t)$ such that  $s_{n}\to s$,
$t_{n}\to t$ as $n\to +\infty$ and
\begin{equation}\label{eq:ynpiumeno1} 
V(u(\cdot,s_{n}))\leq\inf_{y\in (s,s_{n})}V(u(\cdot,y))+\tfrac{1}{n}\hbox{ and }V(u(\cdot,t_{n}))\leq\inf_{y\in (t_{n},t)}V(u(\cdot,y))+\tfrac{1}{n}.
\end{equation}
By (\ref{eq:continuita1}), we have $\|u(\cdot, s_{n})-u(\cdot,s)\|\to 0$ and $\|u(\cdot, t_{n})-u(\cdot,t)\|\to 0$ as $n\to+\infty$, and it is not restrictive to assume that
\begin{equation}\label{eq:ynpiumeno2}
\|u(\cdot, s_{n})-u(\cdot,s)\|\leq d_{0}\hbox{ and }\|u(\cdot, t_{n})-u(\cdot,t)\|\leq d_{0}\hbox{ for any }n\in\N\end{equation}
For every $n\in\N$, consider the paths in $\Gamma$ defined by 
\begin{align*}
 \gamma_{n,-}(y)&=u(\cdot,s)+\frac{y-s}{s_{n}-s}\left(u(\cdot,s_{n})-u(\cdot,s)\right),\quad y\in [s,s_{n}],\\ 
\gamma_{n,+}(y)&=u(\cdot,t)+\frac{t-y}{t-t_{n}}\left(u(\cdot,t_{n})-u(\cdot,t)\right),\quad y\in [t_{n},t].
\end{align*}
Note that, for any $n\in\N$, the paths $\gamma_{n,-}$ and $\gamma_{n,+}$ continuously connect in $\Gamma$ respectively  $u(\cdot,s)$ with $u(\cdot,s_{n})$ and $u(\cdot,t)$ with $u(\cdot,t_{n})$.\\
 Then, since $V(u(\cdot,s)),\, V(u(\cdot,t))\leq c $ and $V(u(\cdot,s_{n})),\, V(u(\cdot,t_{n}))\geq c $, defining for $n\in\N$
\begin{align*}
\sigma_{n}=&\inf\{\bar y\in [s,s_{n}]\, /\, V(\gamma_{n,-}(y))\geq c \hbox{ for every }y\in [\bar y, s_{n}]\},\\
\tau_{n}=&\sup\{\bar y\in [t_{n},t]\, /\, V(\gamma_{n,+}(y))\geq c \hbox{ for every }y\in [t_{n},\bar y]\},
\end{align*}
by continuity, we have that $V(\gamma_{n,-}(\sigma_n))=V(\gamma_{n,+}(\tau_{n}))=c $. Moreover, by definition, 
$V(\gamma_{n,-}(y))\geq c $ for every $y\in [\sigma_n,s_{n}]$ and
$V(\gamma_{n,+}(y))\geq c $ for every $y\in [t_{n},\tau_{n}]$.\hb
Define, for $n,\, j\in\N$,
\begin{equation}\label{eq:w}
w_{n,j}(\cdot,y)=\begin{cases}\gamma_{n,-}(\sigma_n)&\hbox{ if }y\leq \sigma_n,\\
\gamma_{n,-}(y)&\hbox{ if }\sigma_n<y\leq s_{n},\\
u(\cdot,y)&\hbox{ if }s_{n}<y\leq t_{j},\\
\gamma_{j,+}(y)&\hbox{ if }t_{j}<y\leq \tau_{j},\\
\gamma_{j,+}(\tau_{j})&\hbox{ if }\tau_{j}<y,\end{cases}
\end{equation}
and note that, by (\ref{eq:ynpiumeno2}), $w_{n,j}\in\HH_{c}$, and so $\f(w_{n,j})\geq m_{c}$ for all $n,\, j\in\N$.\\
To prove that $\varphi_{c,(s,t)}(u)\geq m_{c}$, since we now know that $\f(w_{n,j})\geq m_{c}$, it will be sufficient to prove that the difference $\varphi_{c,(s,t)}(u)-\f(w_{n,j})$ is definitely nonnegative. To this end, note that, since $\f(w_{n,j})=\varphi_{c,(s,t)}(w_{n,j})$ and since
$w_{n,j}(\cdot,y)=u(\cdot,y)$ for all $y\in (s_{n},t_{j})$, we have
\begin{align*}
\varphi_{c,(s,t)}(u)-\f(w_{n,j})&=\int_{s}^{s_{n}}
\tfrac12\left(\|\partial_{y}u(\cdot,y)\|^{2}-\|\partial_{y}w_{n,j}(\cdot,y)\|^{2}\right)\\ &\qquad
+
\left(V(u(\cdot,y))-V(w_{n,j}(\cdot,y))\right)\, dy\\ &\quad +\int_{t_{j}}^{t}
\tfrac12\left(\|\partial_{y}u(\cdot,y)\|^{2}-\|\partial_{y}w_{n,j}(\cdot,y)\|^{2}\right)\\ &\qquad+
\left(V(u(\cdot,y))-V(w_{n,j}(\cdot,y))\right)\, dy.
\end{align*}
Since $\partial_{y}w_{n,j}(\cdot,y)=\partial_{y}\gamma_{n,-}(\cdot,y)=\frac{1}{s_{n}-s}\left(u(\cdot,s_{n})-u(\cdot,s)\right)$ for all $y\in (\sigma_n,s_{n})$ and $\partial_{y}w_{n,j}(\cdot,y)=0$ for all $y\in (s,s_{{n}})$, by (\ref{eq:continuita1}) we recover that
$$\int_{s}^{s_{n}}
\|\partial_{y}w_{n,j}(\cdot,y)\|^{2}\, dy\leq \frac{1}{s_{n}-s}\|u(\cdot,s_{n})-u(\cdot,s)\|^{2}\leq\int_{s}^{s_{n}}
\|\partial_{y}u(\cdot,y)\|^{2}\, dy.$$
Analogously, we obtain $\int_{t_{j}}^{t}\|\partial_{y}w_{n,j}(\cdot,y)\|^{2}dy\leq\int_{t_{j}}^{t}\|\partial_{y}u(\cdot,y)\|^{2}dy$ and so 
\begin{align}
\varphi_{c,(s,t)}(u)-\f(w_{n,j})&\geq \int_{s}^{s_{n}}V(u(\cdot,y))-V(w_{n,j}(\cdot,y))\, dy \nonumber\\ &\quad+\int_{t_{j}}^{t}V(u(\cdot,y))-V(w_{n,j}(\cdot,y))\, dy.\label{eq:stimadalbasso2}
\end{align}
Then, by (\ref{eq:ynpiumeno1}), Lemma \ref{L:compK} and the continuity property of $V$ it can be proved (see Lemma 3.7 (iv), in \cite{[AlMbrake]}) 
$$
\liminf_{n,j\to+\infty}\,\int_{s}^{s_{n}}V(u(\cdot,y))-V(w_{n,j}(\cdot,y))\, dy  \geq 0
$$ and analogously for the second term $\int_{t_{j}}^{t}V(u(\cdot,y))-V(w_{n,j}(\cdot,y))\, dy$. Hence, by (\ref{eq:stimadalbasso2}) we conclude that $\varphi_{c,(s,t)}(u)=m_c$.

\QED}

Moreover, using Lemmas \ref{L:taglio2(i)}, \ref{L:taglio2(ii)} and \ref{L:taglio2} we have
\begin{Lemma}\label{L:taglio2(iv)}  
There results $\varphi_{c,(s_{c},t_{c})}(u_{c})= m_{c}$ and moreover 
$$
\displaystyle\liminf_{y\to s_{c}^{+}}V(u_{c}(\cdot,y))=\liminf_{y\to t_{c}^{-}}V(u_{c}(\cdot,y))=c .
$$
\end{Lemma}
{\Proof We consider only the case $s_{c},\, t_{c}\in\R$,  similar argument can be used to prove the statement in the cases $s_{c}=-\infty$ or $t_{c}=+\infty$. \\ 
To prove that $\varphi_{c,(s_{c},t_{c})}(u_{c})=m_{c}$, note that by Lemma \ref{L:taglio2(i)}, we have $\varphi_{c,(s_{c},t_{c})}(u_{c})\le m_{c}$. So it will remain to prove that $\varphi_{c,(s_{c},t_{c})}(u_{c})\geq m_{c}$. Note that,  by Lemma \ref{L:taglio2(ii)}, $u(\cdot,t_c)\in\G_c^{+}$, $u(\cdot,s_c)\in\G_{c}^{-} $
and $V(u(\cdot,y))\geq c $ for every $y\in (s_c,t_c)$. Hence, $\varphi_{c,(s_c,t_c)}(u)\geq m_{c}$ follows by Lemma \ref{L:taglio2}.  \\
To prove that $\liminf_{y\to s_{c}^{+}}V(u(\cdot,y))=c $, assume by contradiction that, letting $\ell=\liminf_{y\to s_{c}^{+}}V(u(\cdot,y))$, there results $\ell>c $. Then, considering the sequence $w_{n,j}(\cdot,y)$ defined in (\ref{eq:w}) with $s=s_c$ and $t=t_c$, it can be proved (see Lemma 3.7 (iv), in \cite{[AlMbrake]}) that there exists $\tilde\mu>0$ such that
 for $n$  large and $j\in\N$
\begin{eqnarray*}
\int_{s_{c}}^{s_{n}}V(u_{c}(\cdot,y))-V(w_{n,j}(\cdot,y))\, dy\geq
{\tilde\mu}(s_{n}-s_{c}).
\end{eqnarray*}
Then, by (\ref{eq:stimadalbasso2}) we recover that for $n$ sufficiently large
$$
\varphi_{c,(s_{c},t_{c})}(u_{c})-m_{c}\geq\liminf_{j\to+\infty}\,  [\varphi_{c,(s_{c},t_{c})}(u_{c})-\f(w_{n,j})]\geq {\tilde\mu}(s_{n}-s_{c})>0
$$
which is a contradiction. Analogously, $\liminf_{y\to t_{c}^{-}}V(u(\cdot,y))=c $.\QED}

By Lemma \ref{L:taglio2(iv)} we have that any limit point $u_{c}$  satisfies $\varphi_{c,(s_{c},t_{c})}(u_{c})=m_{c}$ from which we can conclude that every $u_{c}$ is  
a {\sl weak} solution to
(\ref{eq:eq0}) in $\R\times (s_{c},t_{c})$, indeed there result (see e.g. Lemma 3.9 in \cite{[AlMalmost]} or \cite{[AlMbrake]})
        $$
        \displaystyle
    \int_{\R^{2}}\nabla u_{c}\nabla \psi+\nabla W(u_{c})\psi\, dx\, dy=0\quad\hbox{for all $\psi\in C_{0}^\infty(\R\times (s_{c},t_{c}))^{2}$}.
    $$
Then, it is standard to show that $u_{c}$ 
is in fact a {\sl classical} solution to
(\ref{eq:eq0}) on $\R\times (s_{c},t_{c})$.  Hence, considering the {\it energy} 
$$
E_{u}(y)= \tfrac 12\|\partial_{y}u(\cdot,y)\|^{2}_{L^{2}(\R)^{2}}-V(u(\cdot,y))
$$
we get 

\begin{Lemma}\label{L:energy}
There results $E_{u_{c}}(y)=-c $ for all $y\in (s_{c},t_{c})$.
\end{Lemma}
\Proof First note that $E_{u_{c}}(y)$ is constant on $(s_{c},t_{c})$. Indeed, since $u_{c}$ is a { classical} solution to (\ref{eq:eq0}) on $\R\times (s_{c},t_{c})$,  multiplying (\ref{eq:eq0})  by $\partial_{y}u_{c}$ we obtain
\begin{align*}
    0&=
    -\partial_{x,x}u_{c}\cdot \partial_{y}u_{c}-\partial_{y,y}u_{c}\cdot\partial_{y}u_{c}+
    \nabla W (u_{c})\cdot \partial_{y}u_{c}\\
   &= -\partial_{x}(\partial_{x}u_c\cdot
	\partial_{y}u_{c})+\partial_{y}(\tfrac{1}{2}|\partial_{x}u_{c}|^{2}-\tfrac{1}{2}
	|\partial_{y}u_{c}|^{2}+W(u_{c}))\quad\hbox{on }\R\times(s_{c},t_{c}).
\end{align*} 
Then, considering any $(\sigma,\tau)\subset(s_{c},t_{c})$, integrating over $\R\times (\sigma,\tau)$ and using 
Fubini's Theorem, we conclude
\begin{align*}
	   0&=
	-\int_{\sigma}^{\tau}[\int_{\R}\partial_{x}(\partial_{x}u_{c}\cdot
	\partial_{y}u_{c})\,
dx]\,dy +\int_{\R}[\int_{\sigma}^{\tau}\partial_{y}(\tfrac{1}{2}|
	\partial_{x}u_{c}|^{2}-\tfrac{1}{2}|\partial_{y}u_{c}|^{2}+
	W(u_{c}))\, dy]\,dx\\ 
	&=	E_{u_{p}}(\sigma)-E_{u_{p}}(\tau)
\end{align*}
since $\partial_{x}u_{p}(x,y)\cdot
	\partial_{y}u_{p}(x,y)\to 0$ as $x\to\pm\infty$ for almost every $y\in(\sigma,\tau)$.
Moreover we have that for every $(\sigma,\tau)\subset (s_{c}, t_{c})$ there results
\begin{equation}\label{eq:energy2}
\int_{\sigma}^{\tau}\tfrac 12\|\partial_{y}u_{c}(\cdot,y)\|^{2}\, dy=\int_{\sigma}^{\tau} V(u_{c}(\cdot,y))-c \, dy
\end{equation}
Indeed,  
for  $s>0$ let
$$
u_{c}^{s}(\cdot,y)=\begin{cases} u_{c}(\cdot,y+\tau)&y\leq 0,\\
u_{c}(\cdot,\frac{y}{s}+\tau)&y>0\\
\end{cases}
$$
Now note that for every $s>0$, $u_c^s$ verifies the assumption of Lemma \ref{L:taglio2} over $(s_c-\tau, s(t_c-\tau))$, hence 
$\varphi_{c,(s_c-\tau, s(t_c-\tau))}(u^{s}_{c})\geq m_c$ while, by Lemma \ref{L:taglio2(iv)},
$m_c=\varphi_{c,(s_c,t_c)}(u_c)\geq \varphi_{c,(s_c,\tau)}(u_{c})$.
Then we have
\begin{align*}
0\le \varphi_{c,(s_c-\tau, s(t_c-\tau))}(u_{c}^{s})&-\varphi_{c,(s_c,\tau)}(u_{c})
=
(\frac{1}{s}-1)A+
(s-1)B
\end{align*}
where 
$$
A=\displaystyle\int_{\tau}^{t_c}\tfrac{1}{2}\|\partial_{y}u_{c}(\cdot,y)\|^{2}\, dy\quad\hbox{ and}\quad B=\displaystyle\int_{\tau}^{t_c}V(u_{c}(\cdot, y))-c \, dy.
$$
Now note that the real function
$f(s)=(\frac{1}{s}-1)A+(s-1)B$ has minimum value for $s=\sqrt{\frac{A}{B}}$ with
$$
f(\sqrt{\tfrac AB})=-(\sqrt B-\sqrt A)^2
$$
and since $f(s)\geq 0$ for every $s>0$, we conclude $A=B$, that is 
$$
\int_{\tau}^{t_c}\tfrac 12\|\partial_{y}u_{c}(\cdot,y)\|^{2}\, dy=\int_{\tau}^{t_c} V(u_{c}(\cdot,y))-c \, dy.
$$
Analogously  we can prove that 
$$
\int_{s_c}^{\sigma}\tfrac 12\|\partial_{y}u_{c}(\cdot,y)\|^{2}\, dy=\int_{s_c}^{\sigma}V(u_{c}(\cdot,y))-c \, dy
$$
and then, by additivity, we conclude that  (\ref{eq:energy2}) holds.\\
 Then, for every $(\sigma,\tau)\subset(s_{c},t_{c})$  we get 
 $$
 \int_{\sigma}^{\tau}E_{u_{c}}(y)+c \, dy=0
 $$
 and since, as proved above, $E_{u_{c}}(y)$ is constant, the Lemma follows.\QED
 
Now note that by Lemmas \ref{L:taglio2(iv)} and \ref{L:energy}  we plainly obtain that 
  \begin{equation}\label{eq:neumann}
  \liminf_{y\to s_{c}^{+}}\|\partial_{y}u_{c}(\cdot,y)\|= \liminf_{y\to t_{c}^{-}}\|\partial_{y}u_{c}(\cdot,y)\|=0.
    \end{equation}

We will  prove that if $c $ is a regular value of $V$ then
$s_{c},\, t_{c}\in\R$ and in such a case, using the above Neumann conditions, by reflection  we can recover  
an entire periodic solution to (\ref{eq:eq0}), a {\it brake orbit} type solution to (\ref{eq:eq0}). On the other hand, if $c $ is a critical value for $V$ then both $t_{c}$ and $s_{c}$ could be not finite and in such a case from $u_{c}$ we will obtain an {\it homoclinic} or an {\it heteroclinic} type solution to (\ref{eq:eq0}) asymptotic to a critical point of $V$ at level $c $ as $y\to\pm\infty$. First, setting $\KK_{c}^{\pm}=\{q\in\G^{\pm}\,|\, V'(q)=0\hbox{ and }V(q)=c \}$, we have

\begin{Lemma}\label{L:t=infty}
If $t_{c}=+\infty$ then, there exists $q_{0}\in\KK_c^{+}$ such that 
$$
\liminf_{y\to+\infty}\|u_{c}(\cdot,y)-q_{0}\|_{H^{1}}=0.
$$
\end{Lemma}
\Proof
\noindent { We give only a sketch of the proof, we refer to the proof of Lemma 3.11, \cite{[AlMbrake]}, for the details. First, note  that, since $t_{c}=+\infty$, as noted above $u_{c}\in C^{2}(\R\times (s_{c},+\infty))$ is a classical solution to 
(\ref{eq:eq0}) on the half plane $\R\times (s_{c},+\infty)$. Hence, since $\|u_{c}\|_{L^{\infty}(\R^{2})^{2}}\le 1$, using local Schauder estimates, we get  $\norm{u_{c}}_{C^{2}(\R\times (s_{c}+1,+\infty))^{2}}<+\infty$. Fixed any sequence $y_{n}\to +\infty$ and setting $u_{n}(x,y)=u_{c}(x,y+y_{n})$ we 
will prove that, up to a subsequence, $u_{c}(\cdot,y_{n})-q_{0}=u_{n}(\cdot,0)-q_{0}\to 0$ in $H^{1}(\R)^{2}$ for some $q_{0}\in\KK_c^{+}$.\\
\noindent First note that, by regularity, using Lemma \ref{L:taglio2(iv)} and (\ref{eq:neumann}), we obtain that 
\begin{equation}\label{eq:stepN}
\hbox{$\norm{\partial_{y}u_{c}(\cdot,y)}\to 0$ and
$V(u_{c}(\cdot,y))\to c $ as ${y\to +\infty}$.}
\end{equation} 
Using the Ascoli Arzel\`a Theorem we have that there exists $u_{0}\in C^{1}(\R^{2})^{2}$ such that, up to a subsequence, we have $u_{n}\to u_{0}$ in $C^{1}_{loc}(\R^{2})^{2}$. Then, by regularity and (\ref{eq:stepN}) we get that $\partial_{y}u_{0}\equiv 0$ and hence that $u_{0}(x,y)=q_{0}(x)$ for all $(x,y)\in\R^{2}$. Furthermore, since by (\ref{eq:stepN}) we know that  $V(u_{c}(\cdot,y))\to c $ as $y\to+\infty$, using Lemma \ref{L:cattura} we derive that $u_{c}(\cdot,y)\in\G^{+}$ definitively as $y\to+\infty$. In particular $u_{n}(\cdot,0)=u_{c}(\cdot,y_{n})\in\G^{+}$ for $n$ large and by the $C^{1}_{loc}(\R^{2})$ convergence, since $\G^{+}$ is weakly closed, we obtain
 $q_{0}\in\G^{+}$.\\
Finally, since for every $n\in\N$ we have $-\Delta u_{n}+\nabla W(u_{n})=0$ and $u_{n}- q_{0}\to 0$ weakly in $H^{2}(S_{(-1,1)})^{2}$ and in $\CC^{1}_{loc}(\R^{2})^{2}$, by (\ref{eq:stepN}) there results
$$
-\ddot q_{0}(x)+\nabla W(q_{0}(x))=0\hbox{ for every $x\in\R$}
$$
and hence that $V'(q)=0$.

\noindent Then, by Lemma \ref{L:conc}, using (\ref{eq:iper}) and (\ref{eq:stepN}), we obtain that for every $T>0$ we have $u_{n}- q_{0}\to 0$ strongly in $H^{1}(S_{(-T,T)})$.

\noindent In fact, by the Schauder estimates we have that $u_{n}(\cdot,0)-q_{0}\to 0$ in $H^{1}(\R)$ and since $V(u_{n}(\cdot,0))\to c $ by continuity we conclude $V(q_{0})=c $.
\QED}

Analogously we can prove

\begin{Lemma}\label{L:s=infty}
If $s_{c}=-\infty$ then there exists $q_{0}\in\KK_c^{-}$ such that 
$$
\liminf_{y\to-\infty}\|u_{c}(\cdot,y)-q_{0}\|_{H^{1}}=0.
$$ 
\end{Lemma}

By the previous results in particular we obtain that if $t_{c}=+\infty$ or $s_{c}=-\infty$, then $c $ is a critical value for $V$ and hence we obtain

\begin{Corollary}\label{L:ymenoypiufiniti}
If $c $ is a regular value of $V$ then
$s_{c},\, t_{c}\in\R$.
\end{Corollary}

Now, if $s_{c}=-\infty$ and $t_{c}=+\infty$, then the limit point $u_{c}$ is an entire solution to (\ref{eq:eq0}) such that, by Lemmas \ref{L:t=infty} and \ref{L:s=infty}, 
$$
\lim_{y\to\pm\infty}\dd_{H^{1}}(u_{c}(\cdot,y),\KK_{c}^{\pm})=0
$$
that is, $u_{c}$ is an entire solution to (\ref{eq:eq0}) of {\it heteroclinic type}.  As noted in Remark \ref{R:tpm} this is the case that occurs if $c =m$, proving Theorem \ref{T:mainPRIMO}.\smallskip

\noindent If otherwise we have that $s_{c}=-\infty$ and $t_{c}\in\R$, let us consider the function
$$
v_{c}(x,y)=\begin{cases}u_{c}(x,y)&\hbox{if }x\in\R\hbox{ and }y\le t_{c}\\
u_{c}(x,2t_{c}-y)&\hbox{if }x\in\R\hbox{ and }y>t_{c}\end{cases}
$$
Then we have

\begin{Proposition}\label{L:omoclic}
If $s_{c}=-\infty$ and $t_{c}\in\R$, then $v_{c}\in \CC^{2}(\R^{2})$ is a solution of problem (\ref{eq:eq0}). Moreover, $v_{c}(\cdot,t_{c})\in\G_c^{+}$, 
$\partial_{y} v_{c}(\cdot,t_{c})\equiv0$ and there exists $q_{0}\in\KK^{-}_{c}$ such that $\displaystyle\liminf_{y\to\pm\infty}\|v_{c}-q_{0}\|_{H^{1}}=0$.
\end{Proposition}
\Proof By (\ref{eq:neumann}),
there exist four sequences $(\tau^{\pm}_{n})$ and $(\sigma^{\pm}_{n})$ such that
$\sigma^{-}_{n}<\tau^{-}_{n}<t_{c}<\tau^{+}_{n}<\sigma^{+}_{n}$,  for all $n\in\N$, $\tau^{\pm}_{n}\to t_{c}$ and $\sigma^{\pm}_{n}\to \pm\infty$ and 
\begin{equation}\label{eq:sceltatausigma}
\lim_{n\to+\infty}\|\partial_{y}v_{c}(\cdot,\tau^{\pm}_{n})\|= \lim_{n\to+\infty}\|\partial_{y}v_{c}(\cdot,\sigma^{\pm}_{n})\|=0.
\end{equation}
Fixed $\psi\in C_{0}^\infty(\R^{2})$, since $v_{c}$ is a solution to (\ref{eq:eq0}) on $(-\infty,t_{c})$ and $(t_{c},+\infty)$ for  $n$ sufficiently large, by Green's Theorem we obtain
\begin{align*}
    0 = \int_{\R}\int_{\sigma_{n}^{-}}^{\tau_{n}^{-}}&-\Delta v_{c}\,\psi+\nabla W(v^{-}_{c})\psi\,
    dy\, dx\\
    &= \int_{\R}\int_{\sigma_{n}^{-}}^{\tau_{n}^{-}}\nabla v_{c}\nabla\psi +\nabla W(v^{-}_{c})\psi\,
    dy\, dx+\\&\qquad+\int_{\R}\partial_{y}v_{c}(x,\tau_{n}^{-})\psi(x,\tau_{n}^{-})\, dx-\int_{\R}\partial_{y}v_{c}(x,\sigma_{n}^{-})\psi(x,\sigma^{-}_{n})\, dx
    \end{align*}
        and
   \begin{align*}
    0 = \int_{\R}\int_{\tau_{n}^{+}}^{\sigma_{n}^{+}}&-\Delta v_{c}\,\psi+\nabla W(v_{c})\psi\,
    dy\, dx\\
   & = \int_{\R}\int_{\tau_{n}^{+}}^{\sigma_{n}^{+}}\nabla v_{c}\nabla\psi +\nabla W(v_{c})\psi\,
    dy\, dx+\\&\qquad+\int_{\R}\partial_{y}v_{c}(x,\sigma_{n}^{+})\psi(x,\sigma_{n}^{+})\, dx-\int_{\R}\partial_{y}v_{c}(x,\tau_{n}^{+})\psi(x,\tau^{+}_{n})\, dx
    \end{align*}
By (\ref{eq:sceltatausigma}), in the limit for $n\to+\infty$, we obtain that 
$$
\int_{\R}\int_{t_p}^{+\infty}\nabla v_{c}\nabla\psi +\nabla W(v_{c})\psi\,
    dy\, dx=\int_{\R}\int_{-\infty}^{t_p}\nabla v_{c}\nabla\psi +\nabla W(v_{c})\psi\,
    dy\, dx=0$$
Then, $v_{c}$ satisfies 
$$
\int_{\R^{2}}\nabla v_{c}\nabla \psi+\nabla W( v_{c})
     \psi\, dx\, dy=0,\qquad
     \fa \psi\in C_{0}^\infty(\R^{2})
$$
and, hence
we recover that $v_{c}$ 
is in fact a {\sl classical} entire solution to
(\ref{eq:eq0}). \\
Furthermore,  we have that $v_{c}-z_{0}\in H^{2}(S_{(\zeta_{1},\zeta_{2})})$ 
 for every interval $(\zeta_{1},\zeta_{2})\subset \R$ and hence there exists a constant $C>0$ depending only on
$\zeta_{2}-\zeta_{1}$ such that $\norm{v_{c}-z_{0}}_{H^{2}(S_{(\zeta_{1},\zeta_{2})})}\leq C
$.
This implies in particular that
the functions
$y\in \R\to\partial_{y}v_{c}(\cdot,y)\in L^{2}(\R)$ and
$y\in \R\to v_{c}(\cdot,y)\in\Gamma$ are uniformly 
continuous. Then, by Lemma \ref{L:taglio2(ii)}, we have $v_{c}(\cdot,t_{c})\in\G_c^{+}$ and, by Lemma \ref{L:taglio2(iv)},
$$
\|\partial_{y}v_{c}(\cdot,t_{c})\|=\lim_{y\to t_{c}^{-}}\|\partial_{y}v_{c}(\cdot,y)\|=\liminf_{y\to t_{c}^{-}}\|\partial_{y}v_{c}(\cdot,y)\|=
\liminf_{y\to t_{c}^{-}}\|\partial_{y}u_{c}(\cdot,y)\|=0
$$ 
and hence, by continuity, we derive that $\partial_{y} v_{c}(\cdot,t_{c})\equiv0$.
Finally, the asymptotic behaviour of $v_{c}(\cdot, y)$ as $y\to \pm\infty$ follows by Lemma \ref{L:s=infty}.
\QED

Analogously, if $s_{c}\in\R$ and $t_{c}=+\infty$, considering the function
$$
v_{c}((x,y)=\begin{cases}u_{c}(x,y)&\hbox{if }x\in\R\hbox{ and }y\geq s_{c}\\
u_{c}(x,2s_{c}-y)&\hbox{if }x\in\R\hbox{ and }y<s_{c}\end{cases}
$$
we have

\begin{Proposition}\label{L:omoclic2}
If $s_{c}\in\R$ and $t_{c}=+\infty$, then $v_{c}\in \CC^{2}(\R^{2})$ is a solution of problem (\ref{eq:eq0}). Moreover, $v_{c}(\cdot,s_{c})\in\G_c^{-}$, 
$\partial_{y} v_{c}(\cdot,s_{c})\equiv0$ and there exists $q_{0}\in\KK_{c}^{+}$ such that $\displaystyle\liminf_{y\to\pm\infty}\|v_{c}-q_{0}\|_{H^{1}}=0$.
\end{Proposition}

By the previous results, if $s_{c}=-\infty$ and $t_{c}\in\R$ or, respectively, if $s_{c}\in\R$ and $t_{c}=+\infty$ then the corresponding funcion $v_{c}$ is an entire solutions to (\ref{eq:eq0}) of {\it homoclinic type}.\medskip

Finally, if 
$s_{c}, t_{c}\in\R$,  we can define, by reflection and periodic conti\-nua\-tion, a function $v_{c}\in\HH$, periodic in the variable $y$, which we will show to be
an entire solution to (\ref{eq:eq0}), a {\it brake orbits type} solution.\\ Precisely, setting $T_{c}=t_{c}-s_{c}$, let
$$
v_{c}(x,y)=\begin{cases}u_{c}(x,y+s_{c})&\hbox{if }x\in\R\hbox{ and }y\in [0, T_{c})\\
u_{c}(x,t_{c}+T_{c}-y)&\hbox{if }x\in\R\hbox{ and }y\in [T_{c}, 2T_{c}]\end{cases}
$$
and 
$$
v_{c}(x,y)=v_{c}(x,y+2kT_{c})\quad\hbox{ for every }(x,y)\in\R^{2},\ k\in\Z.
$$
Then we have 

\begin{Proposition} \label{L:fin}
If $s_{c}, t_{c}\in\R$, then the function $v_{c}\in C^{2}(\R^{2})$ is a solution of problem (\ref{eq:eq0}). Moreover,
$\partial_{y} v_{c}(\cdot,0)\equiv\partial_{y} v_{c}(\cdot,T_{c})\equiv 0$, $v_{c}(\cdot,0)\in\G_c^{-}$ and $v_{c}(\cdot,T_{c})\in\G_c^{+}$.
\end{Proposition}
\Proof To prove that $v_{c}$ is a weak, and then a classical solution to 
(\ref{eq:eq0}) on $\R$ we can argue as in the proof of Proposition \ref{L:omoclic}, using (\ref{eq:neumann}) and the Green's Formula. By regularity, we have moreover that 
$$
\lim_{y\to 0^{+}}\|\partial_{y}v_{c}(\cdot,y)\|=\liminf_{y\to 0^{+}}\|\partial_{y}v_{c}(\cdot,y)\|=
\liminf_{y\to s_{c}^{+}}\|\partial_{y}u_{c}(\cdot,y)\|=0
$$ 
and analogously, $\lim_{y\to T_{c}^{-}}\|\partial_{y}v_{c}(\cdot,y)\|=0$.
Hence, by continuity we derive that $\partial_{y} v_{c}(\cdot,0)\equiv\partial_{y} v_{c}(\cdot,T_{c})\equiv 0$. Finally, by Lemma \ref{L:taglio2(ii)}, $v_{c}(\cdot,0)\in\G_c^{-}$ and $v_{c}(\cdot,T_{c})\in\G_c^{+}$. \QED

Then,  collecting Propositions \ref{L:omoclic}, \ref{L:omoclic2} and \ref{L:fin}, Theorem  \ref{T:mainPRIMO} and \ref{T:main} follow.

\end{document}